\documentclass[12pt, a4paper]{amsart}
\usepackage{amsmath,amsthm,amsfonts,amscd,amssymb,eucal,latexsym,mathrsfs}
\usepackage[all]{xy}
\usepackage{latexsym, amssymb, amscd, mathrsfs, graphics, fullpage}
\usepackage[backref=page,colorlinks,bookmarksopen=true]{hyperref}

\def\gG{\mathfrak{G}}
\def\gH{\mathfrak{H}}

\def\gM{\mathfrak{M}}
\def\gN{\mathfrak{N}}
\def\gT{\mathfrak{T}}
\def\ga{\mathfrak{a}}
\def\gb{\mathfrak {b}}

\begin{document}

\title[Measured Quantum Groupoid action] 
{Outer actions of Measured Quantum Groupoids}
\author{Michel Enock}
\address{Institut de Math\'ematiques de Jussieu, Unit\'{e} Mixte Paris 6 / Paris 7 /
CNRS de Recherche 7586 \\175, rue du Chevaleret, Plateau 7E, F-75013 Paris}
 \email{enock@math.jussieu.fr}

\begin{abstract}
Mimicking a recent article of Stefaan Vaes, in which was proved that every locally compact quantum group can act outerly, we prove that we get the same result for measured quantum groupoids, with an appropriate definition of outer actions of measured quantum groupoids. This result is used to show that every measured quantum groupoid can be found from some depth 2 inclusion of von Neumann algebras.
 \end{abstract}

\maketitle
\newpage
\section{Introduction}
\label{intro}
\subsection{}
 In two articles (\cite{Val1}, \cite{Val2}), J.-M. Vallin has introduced two notions (pseudo-multiplicative
unitary, Hopf-bimodule), in order to generalize, up to the groupoid
case, the classical notions of multiplicative unitary \cite{BS} and of Hopf-von Neumann algebras \cite{ES}
which were introduced to describe and explain duality of groups, and leaded to appropriate notions
of quantum groups (\cite{ES}, \cite{W1}, \cite{W2}, \cite{BS}, \cite{MN}, \cite{W3}, \cite{KV1}, \cite{KV2}, \cite{MNW}). 
\\ In another article \cite{EV}, J.-M. Vallin and the author have constructed, from a depth 2 inclusion of
von Neumann algebras $M_0\subset M_1$, with an operator-valued weight $T_1$ verifying a regularity
condition, a pseudo-multiplicative unitary, which leaded to two structures of Hopf bimodules, dual
to each other. Moreover, we have then
constructed an action of one of these structures on the algebra $M_1$ such that $M_0$
is the fixed point subalgebra, the algebra $M_2$ given by the basic construction being then
isomorphic to the crossed-product. We construct on $M_2$ an action of the other structure, which
can be considered as the dual action.
\\  If the inclusion
$M_0\subset M_1$ is irreducible, we recovered quantum groups, as proved and studied in former papers
(\cite{EN}, \cite{E2}).
\\ Therefore, this construction leads to a notion of "quantum groupoid", and a construction of a
duality within "quantum groupoids". 
\subsection{}
In a finite-dimensional setting, this construction can be
mostly simplified, and is studied in \cite{NV1}, \cite{BSz1},
\cite{BSz2}, \cite{Sz},\cite{Val3}, \cite{Val4}, \cite{Val5}, and examples are described. In \cite{NV2}, the link between these "finite quantum
groupoids" and depth 2 inclusions of
$II_1$ factors is given, and in \cite{D} had been proved that any finite-dimensional connected $\mathbb{C}^*$-quantum groupoid can act outerly on the hyperfinite $II_1$ factor. 
\subsection{}
In \cite{E3}, the author studied, in whole generality, the notion of pseudo-multiplicative unitary introduced par J.-M. Vallin in \cite{Val2}; following the startegy given by \cite{BS}, with the help of suitable fixed vectors, he introduced a notion of "measured quantum groupoid of compact type". Then F. Lesieur in \cite{L}, using the notion of Hopf-bimodule introduced in \cite{Val1}, then there exist a left-invariant operator-valued weight and a right-invariant operator-valued weight, mimicking in this wider setting the axioms and the technics of Kustermans and Vaes (\cite{KV1}, \cite{KV2}), obtained a pseudo-multiplicative unitary, which, as in the quantum group case, "contains" all the informations about the object (the von Neuman algebra, the coproduct, the antipod, the co-inverse). Lesieur gave the name of "measured quantum groupoids" to these objects. A new set of axioms for these had been given in an appendix of \cite{E5}. In \cite{E4} had been shown that, with suitable conditions, the objects constructed in \cite{EV} from depth 2 inclusions, are "measured quantum groupoids" in the sense of Lesieur.   
\subsection{}
In \cite{E5} have been developped the notions of action (already introduced in \cite{EV}), crossed-product, etc, following what had been done for locally compact quantum groups in (\cite{E1}, \cite{ES1}, \cite{V1}); a biduality theorem for actions had been obtained in (\cite{E5}, 11.6). Moreover, we proved in (\cite{E5} 13.9) that, for any action of a measured quantum groupoid, the inclusion of the initial algebra (on which the measured quantum groupoid is acting) into the crossed-product is depth 2, which leads, thanks to \cite{E4}, to the construction of another measured quantum groupoid (\cite{E5} 14.2).  In \cite{E7} was proved a generalization of Vaes' theorem (\cite{V1}, 4.4) on the standard implementation of an action of a locally compact quantum group; namely, we had obtained such a result when there exists a normal semi-finite faithful operator-valued weight from the von Neumann algebra on which the measured quantum groupoid is acting, onto the copy of the basis of this measured quantum groupoid which is put inside this algebra by the action. 
\subsection{}
One question remained open : can any measured quantum groupoid be constructed from a depth 2 inclusion ? For locally compact quantum groups, the answer is positive (\cite{E5} 14.9), but the most important step in that proof was Vaes' \cite{V2}, who proved that any locally compact quantum group has an outer action. 
\subsection{}
In that article, we answer positively to that question, and we follow the same strategy than for locally compact groups : thanks to the construction of the measured quantum groupoid associated to an action (\cite{E5}, 14.2), we show that this question is equivalent to the existence of an outer action; to prove that last result, we mimick what was done in \cite {V2}, by proving that any measured quantum groupoid has a faithful action. In \cite{V2}, it was constructed on some free product of factors; here we clearly need to construct this action on an amalgated free product of von Neumann algebras, as described, for instance, by Ueda \cite{U}. 
\subsection{}
This article is organized as follows : 
\newline
In chapter \ref{pre}, we recall very quickly all the notations and results needed in that article; we have tried to make these preliminaries as short as possible, and we emphazise that this article should be understood as the continuation of \cite{E5} and \cite{E7}. 
\newline
In chapter \ref{out}, we define outer actions of a measured quantum groupoid, and, in chapter \ref{faithful}, faithful actions and minimal actions, and obtain links between faithful and outer actions. 
\newline
In chapter \ref{anyout}, we construct an outer action of any measured quantum groupoid on some amalgamated free product of von Neumann algebras. 
\newline
Finally, in chapter \ref{S}, we study if and when it is possible for a measured quantum groupoid to act outerly on a semi-finite (or finite) von Neumann algebra, or a finite factor.

 \section{Preliminaries}
 \label{pre}
This article is the continuation of \cite{E5} and \cite{E7}; preliminaries are to be found in \cite{E5}, and we just recall herafter the following definitions and notations :

\subsection{Spatial theory; relative tensor products of Hilbert spaces and fiber products of von Neumann algebras (\cite{C1}, \cite{S}, \cite{T}, \cite{EV})}
\label{spatial}
 Let $N$ be a von Neumann algebra, $\psi$ a normal semi-finite faithful weight on $N$; we shall denote by $H_\psi$, $\gN_\psi$, ... the canonical objects of the Tomita-Takesaki theory associated to the weight $\psi$; let $\alpha$ be a non degenerate faithful representation of $N$ on a Hilbert space $\mathcal H$; the set of $\psi$-bounded elements of the left-module $_\alpha\mathcal H$ is :
\[D(_\alpha\mathcal{H}, \psi)= \lbrace \xi \in \mathcal{H};\exists C < \infty ,\| \alpha (y) \xi\|
\leq C \| \Lambda_{\psi}(y)\|,\forall y\in \gN_{\psi}\rbrace\]
Then, for any $\xi$ in $D(_\alpha\mathcal{H}, \psi)$, there exists a bounded operator
$R^{\alpha,\psi}(\xi)$ from $H_\psi$ to $\mathcal{H}$,  defined, for all $y$ in $\gN_\psi$ by :
\[R^{\alpha,\psi}(\xi)\Lambda_\psi (y) = \alpha (y)\xi\]
which intertwines the actions of $N$. 
\newline
If $\xi$, $\eta$ are bounded vectors, we define the operator product 
\[\langle\xi,\eta\rangle_{\alpha,\psi} = R^{\alpha,\psi}(\eta)^* R^{\alpha,\psi}(\xi)\]
belongs to $\pi_{\psi}(N)'$, which, thanks to Tomita-Takesaki theory, will be identified to the opposite von Neumann algebra $N^o$, which will be equiped with a canonical normal semi-finite faithful weight $\psi^o$. 
\newline
If $y$ in $N$ is analytical with respect to $\psi$, and if $\xi\in D(_\alpha\mathcal H, \psi)$, then we get that $\alpha(y)\xi$ belongs to $D(_\alpha\mathcal H, \psi)$ and that :
\[R^{\alpha, \psi}(\alpha(y)\xi)=R^{\alpha, \psi}(\xi)J_\psi\sigma^\psi_{-i/2}(y^*)J_\psi\]
\newline
If now $\beta$ is a non degenerate faithful antirepresentation of $N$ on a Hilbert space $\mathcal K$, the relative tensor product $\mathcal K\underset{\psi}{_\beta\otimes_\alpha}\mathcal H$ is the completion of the algebraic tensor product $K\odot D(_\alpha\mathcal{H}, \psi)$ by the scalar product defined,  if $\xi_1$, $\xi_2$ are in $\mathcal{K}$, $\eta_1$, $\eta_2$ are in $D(_\alpha\mathcal{H},\psi)$, by the following formula :
\[(\xi_1\odot\eta_1 |\xi_2\odot\eta_2 )= (\beta(\langle\eta_1, \eta_2\rangle_{\alpha,\psi})\xi_1 |\xi_2)\]
If $\xi\in \mathcal{K}$, $\eta\in D(_\alpha\mathcal{H},\psi)$, we shall denote $\xi\underset{\psi}{_\beta\otimes_\alpha}\eta$ the image of $\xi\odot\eta$ into $\mathcal K\underset{\psi}{_\beta\otimes_\alpha}\mathcal H$, and, writing $\rho^{\beta, \alpha}_\eta(\xi)=\xi\underset{\psi}{_\beta\otimes_\alpha}\eta$, we get a bounded linear operator from $\mathcal H$ into $\mathcal K\underset{\nu}{_\beta\otimes_\alpha}\mathcal H$, which is equal to $1_\mathcal K\otimes_\psi R^{\alpha, \psi}(\eta)$. 
\newline
Changing the weight $\psi$ will give a canonical isomorphic Hilbert space, but the isomorphism will not exchange elementary tensors !
\newline

We shall denote $\sigma_\psi$ the relative flip, which is a unitary sending $\mathcal{K}\underset{\psi}{_\beta\otimes_\alpha}\mathcal{H}$ onto $\mathcal{H}\underset{\psi^o}{_\alpha\otimes _\beta}\mathcal{K}$, defined, for any $\xi$ in $D(\mathcal {K}_\beta ,\psi^o )$, $\eta$ in $D(_\alpha \mathcal {H},\psi)$, by :
\[\sigma_\psi (\xi\underset{\psi}{_\beta\otimes_\alpha}\eta)=\eta\underset{\psi^o}{_\alpha\otimes_\beta}\xi\]
In $x\in \beta(N)'$, $y\in \alpha(N)'$, it is possible to define an operator $x\underset{\psi}{_\beta\otimes_\alpha}y$ on $\mathcal K\underset{\psi}{_\beta\otimes_\alpha}\mathcal H$, with natural values on the elementary tensors. As this operator does not depend upon the weight $\psi$, it will be denoted $x\underset{N}{_\beta\otimes_\alpha}y$. We can define a relative flip $\varsigma_N$ at the level of operators such that $\varsigma_N(x\underset{N}{_\beta\otimes_\alpha}y)=y\underset{N^o}{_\alpha\otimes_\beta}x$. If $P$ is a von Neumann algebra on $\mathcal H$, with $\alpha(N)\subset P$, and $Q$ a von Neumann algebra on $\mathcal K$, with $\beta(N)\subset Q$, then we define the fiber product $Q\underset{N}{_\beta*_\alpha}P$ as $\{x\underset{N}{_\beta\otimes_\alpha}y, x\in Q', y\in P'\}'$. 
\newline
It is straightforward to verify that, if $Q_1$ and $P_1$ are two other von Neumann
algebras satisfying the same relations with $N$, we have :
\[Q\underset{N}{_\beta*_\alpha}P\cap Q_1\underset{N}{_\beta*_\alpha}P_1=(Q\cap Q_1)\underset{N}{_\beta*_\alpha}(P\cap P_1)\]
In particular, we have :
\[Q\underset{N}{_\beta*_\alpha} \alpha(N)=(Q\cap\beta (N)')\underset{N}{_\beta\otimes_\gamma} 1\]
\newline
Moreover, this von Neumann algebra can be defined independantly of the Hilbert spaces on which $P$ and $Q$ are represented; if $(i=1,2)$, $\alpha_i$ is a faithful non degenerate homomorphism from $N$ into $P_i$, $\beta_i$ is a faithful non degenerate antihomomorphism from $N$ into $Q_i$, and $\Phi$ (resp. $\Psi$) an homomorphism from $P_1$ to $P_2$ (resp. from $Q_1$ to $Q_2$) such that $\Phi\circ\alpha_1=\alpha_2$ (resp. $\Psi\circ\beta_1=\beta_2$), then, it is possible to define an homomorphism $\Psi\underset{N}{_{\beta_1}*_{\alpha_1}}\Phi$ from $Q_1\underset{N}{_{\beta_1}*_{\alpha_1}}P_1$ into $Q_2\underset{N}{_{\beta_2}*_{\alpha_2}}P_2$. Slice maps with vector states, normal faithful semi-finite weights and operator-valued weights had been defined in \cite{E4} and recalled in (\cite{E5}, 2.5).
\newline
The operators $\theta^{\alpha, \psi}(\xi, \eta)=R^{\alpha, \psi}(\xi)R^{\alpha, \psi}(\eta)^*$, for all $\xi$, $\eta$ in $D(_\alpha\mathcal H, \psi)$, generates a weakly dense ideal in $\alpha(N)'$. Moreover, there exists a family $(e_i)_{i\in I}$ of vectors in $D(_\alpha\mathcal H, \psi)$ such that the operators $\theta^{\alpha, \psi}(e_i, e_i)$ are 2 by 2 orthogonal projections ($\theta^{\alpha, \psi}(e_i, e_i)$ being then the projection on the closure of $\alpha(N)e_i$). Such a family is called an orthogonal $(\alpha, \psi)$-basis of $\mathcal H$.

\subsection{Measured quantum groupoids (\cite{L}, \cite{E5})}
\label{MQG}
 A quintuplet $(N, M, \alpha, \beta, \Gamma)$ will be called a Hopf-bimodule, following (\cite{Val2}, \cite{EV} 6.5), if
$N$,
$M$ are von Neumann algebras, $\alpha$ a faithful non-degenerate representation of $N$ into $M$, $\beta$ a
faithful non-degenerate anti-representation of
$N$ into $M$, with commuting ranges, and $\Gamma$ an injective involutive homomorphism from $M$
into
$M\underset{N}{_\beta *_\alpha}M$ such that, for all $X$ in $N$ :
\newline
(i) $\Gamma (\beta(X))=1\underset{N}{_\beta\otimes_\alpha}\beta(X)$
\newline
(ii) $\Gamma (\alpha(X))=\alpha(X)\underset{N}{_\beta\otimes_\alpha}1$ 
\newline
(iii) $\Gamma$ satisfies the co-associativity relation :
\[(\Gamma \underset{N}{_\beta *_\alpha}id)\Gamma =(id \underset{N}{_\beta *_\alpha}\Gamma)\Gamma\]
This last formula makes sense, thanks to the two preceeding ones and
\ref{spatial}. The von Neumann algebra $N$ will be called the basis of $(N, M, \alpha, \beta, \Gamma)$\vspace{5mm}.\newline
If $(N, M, \alpha, \beta, \Gamma)$ is a Hopf-bimodule, it is clear that
$(N^o, M, \beta, \alpha,
\varsigma_N\circ\Gamma)$ is another Hopf-bimodule, we shall call the symmetrized of the first
one. (Recall that $\varsigma_N\circ\Gamma$ is a homomorphism from $M$ to
$M\underset{N^o}{_r*_s}M$).
\newline
If $N$ is abelian, $\alpha=\beta$, $\Gamma=\varsigma_N\circ\Gamma$, then the quadruplet $(N, M, \alpha, \alpha,
\Gamma)$ is equal to its symmetrized Hopf-bimodule, and we shall say that it is a symmetric
Hopf-bimodule\vspace{5mm}.\newline
A measured quantum groupoid is an octuplet $\mathfrak {G}=(N, M, \alpha, \beta, \Gamma, T, T', \nu)$ such that (\cite{E5}, 3.8) :
\newline
(i) $(N, M, \alpha, \beta, \Gamma)$ is a Hopf-bimodule, 
\newline
(ii) $T$ is a left-invariant normal, semi-finite, faithful operator valued weight $T$ from $M$ to $\alpha (N)$, which means that, for any $x\in\gM_T^+$, we have $(id\underset{\nu}{_\beta*_\alpha}T)\Gamma(x)=T(x)\underset{N}{_\beta\otimes_\alpha}1$. 
\newline
(iii) $T'$ is a right-invariant normal, semi-finite, faithful operator-valued weight $T'$ from $M$ to $\beta (N)$, which means that, for any $x\in\gM_{T'}^+$, we have $(T'\underset{\nu}{_\beta*_\alpha}id)\Gamma(x)=1\underset{N}{_\beta\otimes_\alpha}T'(x)$. 
\newline
(iv) $\nu$ is normal semi-finite faitfull weight on $N$, which is relatively invariant with respect to $T$ and $T'$, which means that the modular automorphisms groups of the weights $\Phi=\nu\circ\alpha^{-1}\circ T$ and $\Psi=\nu^o\circ\beta^{-1}\circ T'$ commute. 
\newline
We shall write $H=H_\Phi$, $J=J_\Phi$, and, for all $n\in N$, $\hat{\beta}(n)=J\alpha(n^*)J$, $\hat{\alpha}(n)=J\beta(n^*)J$.  The weight $\Phi$ will be called the left-invariant weight on $M$. 
\newline
Examples are explained in \ref{d2} and \ref{exMQG}. 
\newline
Then, $\mathfrak {G}$ can be equipped with a pseudo-multiplicative unitary $W$ which is a unitary from $H\underset{\nu}{_\beta\otimes_\alpha}H$ onto $H\underset{\nu^o}{_\alpha\otimes_{\hat{\beta}}}H$ (\cite{E5}, 3.6), which intertwines $\alpha$, $\hat{\beta}$, $\beta$  in the following way : for all $X\in N$, we have :
\[W(\alpha
(X)\underset{N}{_\beta\otimes_\alpha}1)=
(1\underset{N^o}{_\alpha\otimes_{\hat{\beta}}}\alpha(X))W\]
\[W(1\underset{N}{_\beta\otimes_\alpha}\beta
(X))=(1\underset{N^o}{_\alpha\otimes_{\hat{\beta}}}\beta (X))W\]
\[W(\hat{\beta}(X) \underset{N}{_\beta\otimes_\alpha}1)=
(\hat{\beta}(X)\underset{N^o}{_\alpha\otimes_{\hat{\beta}}}1)W\]
\[W(1\underset{N}{_\beta\otimes_\alpha}\hat{\beta}(X))=
(\beta(X)\underset{N^o}{_\alpha\otimes_{\hat{\beta}}}1)W\]
and the operator $W$ satisfies :
\[(1_\gH\underset{N^o}{_\alpha\otimes_{\hat{\beta}}}W)
(W\underset{N}{_\beta\otimes_\alpha}1_{\gH})
=(W\underset{N^o}{_\alpha\otimes_{\hat{\beta}}}1_{\gH})
\sigma^{2,3}_{\alpha, \beta}(W\underset{N}{_{\hat{\beta}}\otimes_\alpha}1)
(1_{\gH}\underset{N}{_\beta\otimes_\alpha}\sigma_{\nu^o})
(1_{\gH}\underset{N}{_\beta\otimes_\alpha}W)\]
Here, $\sigma^{2,3}_{\alpha, \beta}$
goes from $(H\underset{\nu^o}{_\alpha\otimes_{\hat{\beta}}}H)\underset{\nu}{_\beta\otimes_\alpha}H$ to $(H\underset{\nu}{_\beta\otimes_\alpha}H)\underset{\nu^o}{_\alpha\otimes_{\hat{\beta}}}H$, 
and $1_{\gH}\underset{N}{_\beta\otimes_\alpha}\sigma_{\nu^o}$ goes from $H\underset{\nu}{_\beta\otimes_\alpha}(H\underset{\nu^o}{_\alpha\otimes_{\hat{\beta}}}H)$ to $H\underset{\nu}{_\beta\otimes_\alpha}H\underset{\nu}{_{\hat{\beta}}\otimes_\alpha}H$. 
\newline
All the intertwining properties properties allow us to write such a formula, which will be called the
"pentagonal relation". Moreover, $W$, $M$ and $\Gamma$ are related by the following results :
\newline
(i) $M$ is the weakly closed linear space generated by all operators of the form $(id*\omega_{\xi, \eta})(W)$, where $\xi\in D(_\alpha H, \nu)$, and $\eta\in D(H_{\hat{\beta}}, \nu^o)$.
\newline
(ii) for any $x\in M$, we have $\Gamma(x)=W^*(1\underset{N^o}{_\alpha\otimes_{\hat{\beta}}}x)W$. 
\newline
Moreover, it is also possible to construct many other data, namely a co-inverse $R$, a scaling group $\tau_t$, an antipod $S$, a modulus $\delta$, a scaling operator $\lambda$, a managing operator $P$, and a canonical one-parameter group $\gamma_t$ of automorphisms on the basis $N$ (\cite{E5}, 3.8). Instead of $\mathfrak {G}$, we shall mostly use $(N, M, \alpha, \beta, \Gamma, T, RTR, \nu)$ which is another measured quantum groupoid, denoted $\underline{\mathfrak {G}}$, which is equipped with the same data ($W$, $R$, ...) as $\gG$. 
\newline
A dual measured quantum group $\widehat{\mathfrak{G}}$, which is denoted $(N, \widehat{M}, \alpha, \hat{\beta}, \widehat{\Gamma}, \widehat{T}, \widehat{R}\widehat{T}\widehat{R}, \nu)$, can be constructed, and we have $\widehat{\widehat{\mathfrak {G}}}=\underline{\mathfrak {G}}$. 
\newline
Canonically associated to $\mathfrak {G}$, can be defined also the opposite measured quantum groupoid $\mathfrak{G}^o=(N^o, M, \beta, \alpha, \varsigma_N\Gamma, RTR, T, \nu^o)$ and the commutant measured quantum groupoid $\mathfrak{G}^c=(N^o, M', \hat{\beta}, \hat{\alpha}, \Gamma^c, T^c, R^cT^cR^c, \nu^o)$; we have $(\mathfrak{G}^o)^o=(\mathfrak{G}^c)^c=\underline{\mathfrak{G}}$, $\widehat{\mathfrak{G}^o}=(\widehat{\mathfrak {G}})^c$, $\widehat{\mathfrak {G}^c}=(\widehat{\mathfrak {G}})^o$, and $\mathfrak{G}^{oc}=\mathfrak {G}^{co}$ is canonically isomorphic to $\underline{\mathfrak {G}}$ (\cite{E5}, 3.12). 
\newline
The pseudo-multiplicative unitary of $\widehat{\mathfrak{G}}$ (resp. $\mathfrak{G}^o$, $\mathfrak{G}^c$) will be denoted $\widehat{W}$ (resp. $W^o$, $W^c$). The left-invariant weight on $\widehat{\mathfrak{G}}$ (resp. $\mathfrak{G}^o$, $\mathfrak{G}^c$) will be denoted $\widehat{\Phi}$ (resp. $\Phi^o$, $\Phi^c$). 
\newline
Let $_a\gH_b$ be a $N-N$-bimodule, i.e. an Hilbert space $\gH$ equipped with a normal faithful non degenerate representation $a$ of $N$ on $\gH$ and a normal faithful non degenerate anti-representation $b$ on $\gH$, such that $b(N)\subset a(N)'$. A corepresentation of $\gG$ on $_a\gH_b$ is a unitary $V$ from $\gH\underset{\nu^o}{_a\otimes_\beta}H$ onto 
 $\gH\underset{\nu}{_b\otimes_\alpha}H$, satisfying, for all $n\in N$ :
 \[V(b(n)\underset{N^o}{_a\otimes_\beta}1)=(1\underset{N}{_b\otimes_\alpha}\beta(n))V\]
 \[V(1\underset{N^o}{_a\otimes_\beta}\alpha(x))=(a(n)\underset{N}{_b\otimes_\alpha}1)V\]
such that, for any $\xi\in D(_a\gH, \nu)$ and $\eta\in D(\gH_b, \nu^o)$, the operator $(\omega_{\xi, \eta}*id)(V)$ belongs to $M$ (then, it is possible to define $(id*\theta)(V)$, for any $\theta$ in $M_*^{\alpha, \beta}$ which is the linear set generated by the $\omega_\xi$, with $\xi\in D(_\alpha H, \nu)\cap D(H_\beta, \nu^o)$), and such that the application $\theta\rightarrow (id*\theta)(V)$ from $M_*^{\alpha, \beta}$ into $\mathcal L(\gH)$ is multiplicative (\cite{E5} 5.1, 5.5).

\subsection{Action of a measured quantum groupoid (\cite{E5})}
\label{action}

An action (\cite{E5}, 6.1) of $\mathfrak{G}$ on a von Neumann algebra $A$ is a couple $(b, \mathfrak a)$, where :
\newline
(i) $b$ is an injective $*$-antihomomorphism from $N$ into $A$; 
\newline
(ii) $\mathfrak a$ is an injective $*$-homomorphism from $A$ into $A\underset{N}{_b*_\alpha}M$; 
\newline
(iii) $b$ and $\mathfrak a$ are such that, for all $n$ in $N$:
\[\mathfrak a (b(n))=1\underset{N}{_b\otimes_\alpha}\beta(n)\]
(which allow us to define $\mathfrak a\underset{N}{_b*_\alpha}id$ from $A\underset{N}{_b*_\alpha}M$ into $A\underset{N}{_b*_\alpha}M\underset{N}{_\beta*_\alpha}M$)
and such that :
\[(\mathfrak a\underset{N}{_b*_\alpha}id)\mathfrak a=(id\underset{N}{_b*_\alpha}\Gamma)\mathfrak a\]
The invariant subalgebra $A^\ga$ is defined by :
\[A^\ga=\{x\in A\cap b(N)'; \ga(x)=x\underset{N}{_b\otimes_\alpha}1\}\]
Let us write, for any $x\in A^+$, $T_\ga(x)=(id\underset{\nu}{_b*_\alpha}\Phi)\ga(x)$; this formula defines a normal faithful operator-valued weight from $A$ onto $A^\ga$; the action $\ga$ will be said integrable if $T_\ga$ is semi-finite (\cite{E5}, 6.11, 12, 13 and 14). 
\newline
If the von Neumann algebra acts on a Hilbert space $\gH$, and if there exists a representation $a$ of $N$ on $\gH$ such that $b(N)\subset A\subset a(N)'$, a corepresentation $V$ of $\gG$ on the bimodule $_a\gH_b$ will be called an implementation of $\ga$ if we have $\ga(x)=V(x\underset{N^o}{_a\otimes_b}1)V^*$ , for all $x\in A$ (\cite{E5}, 6.6).

\subsection{Crossed-product (\cite{E5})}
\label{crossed}
The crossed-product of $A$ by $\mathfrak {G}$ via the action $\mathfrak a$ is the von Neumann algebra generated by $\mathfrak a(A)$ and $1\underset{N}{_b\otimes_\alpha}\widehat{M}'$ (\cite{E5}, 9.1) and is denoted $A\rtimes_\mathfrak a\mathfrak {G}$; then there exists (\cite{E5}, 9.3) an integrable action $(1\underset{N}{_b\otimes_\alpha}\hat{\alpha}, \tilde{\mathfrak a})$ of $(\widehat{\mathfrak {G}})^c$ on $A\rtimes_\mathfrak a\mathfrak {G}$, called the dula action. 
\newline
The biduality theorem (\cite{E5}, 11.6) says that the bicrossed-product $(A\rtimes_\mathfrak a\mathfrak {G})\rtimes_{\tilde{\mathfrak a}}\widehat{\mathfrak {G}}^o$ is canonically isomorphic to $A\underset{N}{_b*_\alpha}\mathcal L(H)$; more precisely, this isomorphism is given by :
\[\Theta (\ga\underset{N}{_b*_\alpha}id)(A\underset{N}{_b*_\alpha}\mathcal L(H))=(A\rtimes_\mathfrak a\mathfrak {G})\rtimes_{\tilde{\mathfrak a}}\widehat{\mathfrak {G}}^o\]
where $\Theta$ is the spatial isomorphism between $\mathcal L(\gH\underset{\nu}{_b\otimes_\alpha}H\underset{\nu}{_\beta\otimes_\alpha}H)$ and $\mathcal L(\gH\underset{\nu}{_b\otimes_\alpha}H\underset{\nu^o}{_{\hat{\alpha}}\otimes_\beta}H)$ implemented by $1_\gH\underset{\nu}{_b\otimes_\alpha}\sigma_\nu W^o\sigma_\nu$; the biduality theorem says also that this isomorphism sends  the action $(1\underset{N}{_b\otimes_\alpha}\hat{\beta}, \underline{\mathfrak a})$ of $\gG$ on $A\underset{N}{_b*_\alpha}\mathcal L(H)$, defined, for any $X\in A\underset{N}{_b*_\alpha}\mathcal L(H)$, by :
\[\underline{\mathfrak a}(X)=(1\underset{N}{_b\otimes_\alpha}\sigma_{\nu^o}W\sigma_{\nu^o})(id\underset{N}{_b*_\alpha}\varsigma_N)(\mathfrak a\underset{N}{_b*_\alpha}id)(X)(1\underset{N}{_b\otimes_\alpha}\sigma_{\nu^o}W\sigma_{\nu^o})^*\]
on the bidual action (of $\mathfrak{G}^{co}$) on $(A\rtimes_\mathfrak a\mathfrak {G})\rtimes_{\tilde{\mathfrak a}}\widehat{\mathfrak {G}}^o$. 
\newline
We have $(A\rtimes_\ga\gG)^{\tilde{\ga}}=\ga(A)$ (\cite{E5} 11.5), and, therefore, the normal faithful semi-finite operator-valued weight $T_{\tilde{\ga}}$ sends $A\rtimes_\mathfrak a\mathfrak {G}$ onto $\ga(A)$; therefore, starting with a  normal semi-finite weight $\psi$ on $A$, we can construct a dual weight $\tilde{\psi}$ on $A\rtimes_\mathfrak a\mathfrak {G}$ by the formula $\tilde{\psi}=\psi\circ\ga^{-1}\circ T_{\tilde{\ga}}$ (\cite{E5} 13.2). 
\newline
Moreover (\cite{E5} 13.3), the linear set generated by all the elements $(1\underset{N}{_b\otimes_\alpha}a)\mathfrak a(x)$, for all $x\in\gN_\psi$, $a\in\gN_{\widehat{\Phi}^c}\cap\gN_{\hat{T}^c}$, is a core for $\Lambda_{\tilde{\psi}}$, and it is possible to identify the GNS representation of $A\rtimes_\mathfrak a\gG$ associated to the weight $\tilde{\psi}$ with the natural representation on $H_\psi\underset{\nu}{_b\otimes_\alpha}H_\Phi$ by writing :
\[\Lambda_\psi(x)\underset{\nu}{_b\otimes_\alpha}\Lambda_{\widehat{\Phi}^c}(a)=\Lambda_{\tilde{\psi}}[(1\underset{N}{_b\otimes_\alpha}a)\mathfrak a(x)]\]
which leads to the identification of $H_{\tilde{\psi}}$ with $H_\psi\underset{\nu}{_b\otimes_\alpha}H$. Let us write, for all $n\in N$, $a(n)=J_\psi b(n^*)J_\psi$.
\newline
Then, the unitary $U_\psi^\ga=J_{\tilde{\psi}}(J_\psi\underset{N^o}{_a\otimes_\beta}J_{\widehat{\Phi}})$ from $H_\psi\underset{\nu^o}{_a\otimes_\beta}H_\Phi$ onto $H_\psi\underset{\nu}{_b\otimes_\alpha}H_\Phi$ satisfies :
\[U^\mathfrak a_\psi(J_\psi\underset{N}{_b\otimes_\alpha}J_{\widehat{\Phi}})=(J_\psi\underset{N}{_b\otimes_\alpha}J_{\widehat{\Phi}})(U^\mathfrak a_\psi)^*\]
and we have (\cite{E5} 13.4) :
\newline
(i) for all $y\in A$ :
\[\mathfrak a (y)=U^\mathfrak a_\psi(y\underset{N^o}{_a\otimes_\beta}1)(U^\mathfrak a_\psi)^*\]
(ii) for all $b\in M$ :
\[(1\underset{N}{_b\otimes_\alpha}J_\Phi bJ_\Phi)U^\mathfrak a_\psi=U^\mathfrak a_\psi(1\underset{N^o}{_a\otimes_\beta}J_\Phi bJ_\Phi)\]
(iii) for all $n\in N$ :
\[U_\psi^\mathfrak a(b(n)\underset{N^o}{_a\otimes_\beta}1)=(1\underset{N}{_b\otimes_\alpha}\beta(n))U_\psi^\mathfrak a\]
\[U_\psi^\mathfrak a(1\underset{N^o}{_a\otimes_\beta}\alpha(n))=(a(n)\underset{N}{_b\otimes_\alpha}1)U_\psi^\mathfrak a\]
Finally, if there exists a normal faithful semi-finite operator-valued weight $\gT$ from $A$ on $b(N)$ such that $\psi=\nu^o\circ b^{-1}\circ \gT$, then, we can prove (\cite{E7} 5.7 and 5.8) that $U_\ga^\psi$ is a corepresentation, which, by (i), implements $\ga$, that we shall call a standard implementation of $\ga$. 
\newline
In \cite{E7} was introduced the notion of invariant weight with respect to an action; namely, be given an action $(b, \ga)$ of a measured quantum groupoid $\gG$ on a von Neumann algebra $A$, a normal faithful semi-finite weight $\psi$ on $A$ is said invariant by $\ga$, if, for any $\eta\in D(_\alpha H, \nu)\cap D(H_\beta, \nu^o)$ and $x\in\gN_\psi$, we have :
\[\psi[(id\underset{N}{_b*_\alpha}\omega_\eta)\ga(x^*x)]=\|\Lambda_\psi(x)\underset{\nu^o}{_a\otimes_\beta}\eta\|^2\]
and if $\psi$ bear a density property (\cite{E5}, 8.1), namely that $D((H_\psi)_b, \nu^o)\cap D(_a H_\psi, \nu)$ is dense in $H_\psi$. 
\newline
Then, it was proved (\cite{E7}, 7.7 (vi)) that there exists a normal semi-finite faithful operator-valued weight $\gT$ from $A$ on $b(N)$, such that $\psi=\nu^o\circ b^{-1}\circ \gT$ (and, therefore, $U_\ga^\psi$ is a standard implementation of $\ga$); moreover, $\gT$ satisfies, for all $x\in\gN_\psi\cap\gN_\gT$ :
\[(\gT\underset{N}{_b*_\alpha}id)\ga(x^*x)=1\underset{N}{_b\otimes_\alpha}\beta\circ b^{-1}\gT(x^*x)=\ga(\gT(x^*x))\]
and $(\psi\underset{\nu}{_b*_\alpha}id)\ga(x^*x)=\beta\circ b^{-1}\gT(x^*x)$. Such an operator-valued weight $\gT$ will be called invariant under $\ga$. A normal faithful conditional expectation $E$ from $A$ to $b(N)$ will be called invariant if $(E\underset{N}{_b*_\alpha}id)\ga=\ga\circ E$. 

\subsection{Depth 2 inclusions}
\label{d2}
Let $M_0\subset M_1$ be an inclusion of $\sigma$-finite von Neumann algebras, equipped with a normal faithful semi-finite operator-valued weight $T_1$ from $M_1$ to $M_0$ (to be more precise, from $M_1^{+}$ to the extended positive elements of $M_0$ (cf. \cite{T} IX.4.12)). Let $\psi_0$ be a normal faithful semi-finite weight on $M_0$, and $\psi_1=\psi_0\circ T_1$; for $i=0,1$, let $H_i=H_{\psi_i}$, $J_i=J_{\psi_i}$, $\Delta_i=\Delta_{\psi_i}$ be the usual objects constructed by the Tomita-Takesaki theory associated to these weights. We shall write $j_i$ for the mirroring on $\mathcal L(H_i)$ defined by $j_i(x)=J_ix^*J_i$. We shall write also $j_1$ for the restriction of the mirroring to $M'_0\cap M_2$ (which is an anti-automorphism of $M'_0\cap M_2$), or for the restriction of the mirroring to $M'_0\cap M_1$ (which is an injective anti-homomorphism from $M'_0\cap M_1$ into $M'_0\cap M_2$).
\newline
 Following (\cite{J}, 3.1.5(i)), the von Neumann algebra $M_2=J_1M'_0J_1$ defined on the Hilbert space $H_1$ will be called the basic construction made from the inclusion $M_0\subset M_1$. We have $M_1\subset M_2$, and we shall say that the inclusion $M_0\subset M_1\subset M_2$ is standard.  
 \newline
 Using Haagerup's theorem (\cite{T}, 4.24), we can construct from $T_1$ another normal faithful semi-finite operator-valued weight $T'_1$ from $M'_0$ onto $M'_1$, and, by defintion of $M_2$, a normal faithful semi-finite operator-valued weight $T_2$ from $M_2$ onto $M_1$; $T_2$ will be called the basic constrction made from $T_1$; we can go on and construct $M_2\subset M_3$ and $T_3$ by the basic construction made from $M_1\subset M_2$ and $T_2$. 
 \newline
Following now (\cite{GHJ} 4.6.4), we shall say that the inclusion $M_0\subset M_1$ is depth 2 if the inclusion (called the derived tower) :
\[M'_0\cap M_1\subset M'_0\cap M_2\subset M'_0\cap M_3\]
 is also standard, and, following (\cite{EN}, 11.12), we shall say that the operator-valued weight $T_1$ is regular if both restrictions $\tilde{T_2}=T_{2|M'_0\cap M_2}$ and $\tilde{T_3}=T_{3|M'_1\cap M_3}$ are semi-finite. 
 \newline
 In \cite{EV} was proved that, with such an hypothesis, there exists a coproduct $\Gamma$ from $M'_0\cap M_2$ into $(M'_0\cap M_2)\underset{M'_0\cap M_1}{_{j_1}*_{id}}(M'_0\cap M_2)$ (where here $id$ means the injection of $M'_0\cap M_1$ into $M'_0\cap M_2$) such that 
 $(M'_0\cap M_1, M'_0\cap M_2, id, j_1, \Gamma)$
 is a Hopf-bimodule; moreover, $\tilde{T_2}$ is then a left-invariant weight, and $j_1\circ \tilde{T_2}\circ j_1$ a right-invariant weight; if there exists a normal faithful semi-finite weight $\chi$ on $M'_0\cap M_1$ invariant under the modular automorphism group $\sigma_t^{T_1}$ (\cite{E4} 8.2 and 8.3), we get that  :
 \[\gG_1=(M'_0\cap M_1, M'_0\cap M_2, id, j_1, \Gamma, \tilde{T_2}, j_1\circ \tilde{T_2}\circ j_1, \chi)\]
 is a measured quantum groupoid. We shall write $\gG_1=\gG(M_0\subset M_1)$. 
 \newline
 Moreover, the inclusion $M_1\subset M_2$ satisfies the same hypothesis, and leads to another measured quantum groupoid $\gG_2$, which can be identified with $\widehat{\gG_1}^o$, and there exists a canonical action $\ga$ of $\gG_2$ on $M_1$ (\cite{EV}, 7.3), which can be described as follows : the anti-representation of the basis $M'_1\cap M_2$ (which, using $j_1$, is anti-isomorphic to $M'_0\cap M_1$), is given by the natural inclusion of $M'_0\cap M_1$ into $M_1$, and the homomorphism from $M_1$ is given by the natural inclusion of $M_1$ into $M_3$ (which is, thanks to (\cite{EV}, 4.6), isomorphic to $M_1\underset{M'_0\cap M_1}{_{j_1}*_{id}}\mathcal L(H_{\chi_2})$, where $\chi_2=\chi\circ\tilde{T_2}$). We then get that $M_0=M_1^{\ga}$ and that $M_2$ is isomorphic to $M_1\rtimes_\ga\gG_2$ (\cite{EV}, 7.5 and 7.6). 
\newline
So, from a depth 2 inclusion $M_0\subset M_1$ equipped with a regular operator-valed weight, and an invariant weight on the first relative commutant, one can construct a measured quantum groupoid $\gG_2$ given, in fact, by 
a specific action $\ga$ of $\gG_2$ on $M_1$, with $M_0$ being the invariant elements under this action. 
\newline
If $\gG$ is any measured quantum groupoid, and $(b, \ga)$ an action of $\gG$ on a von Neumann algebra $A$; then the inclusion $\ga(A)\subset A\rtimes_\ga \gG$ is depth 2 (\cite{E5}, 13.9), and the operator-valued weight $T_{\tilde{\ga}}$ is regular (\cite{E5}, 13.10); so, we can construct a Hopf-bimodule from this depth 2 inclusion, equipped with a left-invariant operator-valued weight and a right-invariant operator-valued weight; moreover, if there exists a weight $\chi$ on $A\rtimes_\ga \gG\cap \ga(A)'$, invariant by $\sigma_t^{T_{\tilde{\ga}}}$, we get another measured quantum groupoid $\gG(\ga)=\gG(\ga(A)\subset A\rtimes_\ga\gG)$ (\cite{E5}, 14.2), which contains, in a sense, $\gG^{oc}$ (\cite{E5}, 14.7). 
\newline
More precisely, the inclusion $\ga(A)\subset A\rtimes_\ga\gG\subset A\underset{N}{_b*_\alpha}\mathcal L(H)$ is standard, and, if we write $B=A\rtimes_\ga \gG\cap \ga(A)'$ and $\gb=\tilde{\ga}_{|B}$, the derived inclusion $B\subset A\underset{N}{_b*_\alpha}\mathcal L(H)\cap\ga(A)'$ is isomorphic to $\gb(B)\subset B\rtimes_\gb\widehat{\gG}^c$ (\cite{E5}, 13.9), and there exist a $*$-anti-automorhism $j_1$ of $B\rtimes_\gb\widehat{\gG}^c$ and a coproduct $\Gamma_1$ such that (\cite{E5} 14.2) :
\[\gG(\ga)=(B, B\rtimes_\gb\widehat{\gG}^c, \gb, j_1\circ\gb, \Gamma_1, T_{\tilde{\gb}}, j_1\circ T_{\tilde{\gb}}\circ j_1, \chi)\]

\subsection{Examples of measured quantum groupoids}
\label{exMQG}
Examples of measured quantum groupoids are the following :
\newline
(i) locally compact quantum groups, as defined and studied by J. Kustermans and S. Vaes (\cite{KV1}, \cite {KV2}, \cite{V1}); these are, trivially, the measured quantum groupoids with the basis $N=\mathbb{C}$. 
\newline
(ii) measured groupoids, equipped with a left Haar system and a quasi-invariant measure on the set of units, as studied mostly by T. Yamanouchi (\cite{Y1}, \cite{Y2}, \cite{Y3}, \cite{Y4}); it was proved in \cite{E6} that these measured quantum groupoids are exactly those whose underlying von Neumann algebra is abelian. 
\newline
(iii) the finite dimensional case had been studied by D. Nikshych and L. Vainermann (\cite{NV1}, \cite{NV2}, \cite{NV3}), J.-M. Vallin (\cite{Val3}, \cite{Val4}) and M.-C. David (\cite{D}); in that case, non trivial examples are given, for instance Temperley-Lieb algebras (\cite{NV3}, \cite{D}), which had appeared in subfactor theory (\cite{J}).  . 
\newline
(iv) continuous fields of ($\bf{C}^*$-version of) locally compact quantum groups, as studied by E. Blanchard in (\cite{Bl1}, \cite{Bl2}); it was proved in \cite{E6} that these measured quantum groupoids are exactly those whose basis is central in the underlying von Neumann algebras of both the measured quantum groupoid and its dual. As a particular case, we find in (\cite{L}, 17.1) that, be given a family $\gG_i=(N_i, M_i, \alpha_i, \beta_i, \Gamma_i, T_i, T'_i, \nu_i)$ a measured quantum groupoids, Lesieur showed that it is possible to construct another measured quantum groupoid $\gG=\oplus_{i\in I}\gG_i=(\oplus_{i\in I}N_i, \oplus_{i\in I}M_i, \oplus_{i\in I}\alpha_i, \oplus_{i\in I}\beta_i, \oplus_{i\in I}\Gamma_i, \oplus_{i\in I}T_i, \oplus_{i\in I}T'_i, \oplus_{i\in I}\nu_i)$. 
\newline
(v) in \cite{DC}, K. De Commer proved that, in the case of a monoidal equivalence between two locally compact quantum groups (which means that these two locally compact quantum group have commuting ergodic and integrable actions on the same von Neumann algebra), it is possible to construct a measurable quantum groupoid of basis $\mathbb{C}^2$ which contains all the data. Moreover, this construction was usefull to prove new results on locally compact quantum groups, namely on the deformation of a locally compact quantum group by a unitary $2$-cocycle; he proved that these measured quantum groupoids are exactly those whose basis $\mathbb{C}^2$ is central in the underlying von Neumann algebra of the measured quatum groupoid, but not in the underlying von Neumann algebra of the dual measured quantum groupoid. 
\newline
(vi) in \cite{VV} and \cite{BSV} was given a specific procedure for constructing locally compact quantum groups, starting from a locally compact group $G$, whose almost all elements belong to the product $G_1G_2$ (where $G_1$ and $G_2$ are closed subgroups of $G$ whose intersection is reduced to the unit element of $G$); such $(G_1, G_2)$ is called a "matched pair" of locally compact groups (more precisely, in \cite{VV}, the set $G_1G_2$ is required to be open, and it is not the case in \cite{BSV}).Then, $G_1$ acts naturally on $L^\infty(G_2)$ (and vice versa), and the two crossed-products obtained bear the structure of two locally compact quantum groups in duality. In \cite{Val5}, J.-M. Vallin generalizes this constructions up to groupoids, and, then, obtains examples of measured quantum groupoids; more specific examples are then given by the action of a matched pair of groups on a locally compact space, and also more exotic examples.

\section{Outer actions of a measured quantum groupoid}
\label{out}
In this chapter, we define (\ref{defout}) outer actions of a measured quantum groupoid, and prove (\ref{thgeo}) that a measured quantum groupoid can be constructed by a geometric construction from a depth 2 inclusion if and only if it has an outer action on some von Neumann algebra. 

\subsection{Theorem}
\label{thout}
{\it Let $\gG=(N, M, \alpha, \beta, \Gamma, T, T', \nu)$ be a measured quantum group, and $(b, \ga)$ an action of $\gG$ on a von Neumann algebra $A$; then, are equivalent :
\newline
(i) $A\rtimes_\ga\gG\cap \ga(A)'=1\underset{N}{_b\otimes_\alpha}\hat{\alpha}(N)$; 
\newline
(ii) $A\underset{N}{_b*_\alpha}\mathcal L(H)\cap \ga(A)'=1\underset{N}{_b\otimes_\alpha}M'$; 
\newline
(iii) it is possible to define the measured quantum groupoid $\gG(\ga)$, and $\gG(\ga)=\gG^{oc}$. }

\begin{proof}
Let us suppose (i); using then (\cite{E5}, 14.1 (iii)), we get (ii). 
\newline
Let us suppose (ii); using (\cite{E5}, 14.7), we see the application $x\mapsto 1\underset{N}{_b\otimes_\alpha}x$ from $M'$ onto $A\underset{N}{_b*_\alpha}\mathcal L(H)\cap \ga(A)'$ is an isomorphism of Hopf-bimodules, from $\gG^{oc}$ onto the Hopf-bimodule constructed from the depth 2 inclusion $\ga(A)\subset A\rtimes_\ga\gG$, and sends the left- (resp. right-) invariant operator-valued weights of $\gG^{oc}$ on the left(resp. right-) invariant operator-valued weights of Hopf-bimodule constructed from the depth 2 inclusion $\ga(A)\subset A\rtimes_\ga\gG$; therefore, we get that the weight $\nu$ is invariant under $\sigma_t^{T_{\tilde{\ga}}}$, which means that we can define the measured quantum groupoid $\gG(\ga)$, and that $\gG(\ga)=\gG^{oc}$, which is (iii).
\newline
Let us suppose (iii); the application $x\mapsto 1\underset{N}{_b\otimes_\alpha}x$ from $M'$ onto $A\underset{N}{_b*_\alpha}\mathcal L(H)\cap \ga(A)'$ is an isomorphism between $\gG^{oc}$ and $\gG(\ga)$; in particular, we get (i). \end{proof}

\subsection{Definition}
\label{defout}
Let $\gG=(N, M, \alpha, \beta, \Gamma, T, T', \nu)$ be a measured quantum groupoid, and $(b, \ga)$ an action of $\gG$ on a von Neumann algebra $A$; we shall say that the action $(b, \ga)$ is outer if it satisfies one of the equivalent conditions of \ref{thout}. 

\subsection{Theorem}
\label{dualout}
{\it Let $\gG=(N, M, \alpha, \beta, \Gamma, T, T', \nu)$ be a measured quantum groupoid, and $(b, \ga)$ an outer action of $\gG$ on a von Neumann algebra $A$; then, the dual action of the measured quantum groupoid $\widehat{\gG}^c$ on the crossed product $A\rtimes_\ga\gG$ is outer. }

\begin{proof}
Let us put the von Neumann algebra $A$ on its standard Hilbert space $L^2(A)$; we have, using \ref{thout} and (\cite{E5}, 3.11) :
\begin{eqnarray*}
A\underset{N}{_b*_\alpha}\mathcal L (H)\cap (A\rtimes_\ga\gG)'
&=&
A\underset{N}{_b*_\alpha}\mathcal L (H)\cap\ga(A)'\cap \mathcal L (L^2(A))\underset{N}{_b*_\alpha}\widehat{M}\\
&=&
1\underset{N}{_b\otimes_\alpha}M'\cap \mathcal L(L^2(A))\underset{N}{_b*_\alpha}\widehat{M}\\
&=&
1\underset{N}{_b\otimes_\alpha}M'\cap \widehat{M}\\
&=&
1\underset{N}{_b\otimes_\alpha}\hat{\beta}(N)
\end{eqnarray*}
from which we get the result, using \ref{thout} again. \end{proof}

\subsection{Proposition}
\label{propout}
{\it Let $\gG=(N, M, \alpha, \beta, \Gamma, T, T', \nu)$ be a measured quantum group, and $(b, \ga)$ an outer action of $\gG$ on a von Neumann algebra $A$; then, we have :
\[Z(A)=\{b(n), n\in Z(N), \beta(n)\in Z(M)\}\]
Moreover, we have :}
\[Z(A\rtimes_\ga\gG)=\{1\underset{N}{_b\otimes_\alpha}\hat{\alpha}(n), \alpha(n)\in Z(\widehat{M})\}\]

\begin{proof}
As we have $\ga(Z(A))\subset A\rtimes_\ga\gG\cap \ga(A)'$, we get that, for any $z\in Z(A)$, there exists $n\in N$ such that $\ga(z)=1\underset{N}{_b\otimes_\alpha}\hat{\alpha}(n)$. But, we then infer that $\hat{\alpha}(n)$ belongs to $M\cap M'\cap \widehat{M}'$; therefore, we have $\hat{\alpha}(n)=\beta(n)\in Z(M)$, $n\in Z(N)$, and $\ga(z)=1\underset{N}{_b\otimes_\alpha}\beta(n)=\ga(b(n))$, from which we gat that 
\[Z(A)\subset\{b(n), n\in Z(N), \beta(n)\in Z(M)\}\]
Conversely, if $n\in Z(N)$, such that $\beta(n)\in Z(M)$, we get that $\ga(b(n))=1\underset{N}{_b\otimes_\alpha}\beta(n)$ commutes with all elements $\ga(x)\in A\underset{N}{_b*_\alpha}M$, for any $x\in A$; therefore, we get that $b(n)\in Z(A)$. Applying this result to the outer action $\tilde{\ga}$, we get that :
\[Z(A\rtimes_\ga\gG)=\{1\underset{N}{_b\otimes_\alpha}\hat{\alpha}(n), \hat{\alpha}(n)\in Z(\widehat{M}')\}\]
and, as $\hat{\alpha}(n)=J\hat{J}\alpha(n)\hat{J}J$, where $\hat{J}$ stands for $J_{\widehat{\Phi}}$, (\cite{E5}, 3.11), we get the result. 
\end{proof}

\subsection{Corollary and Definition}
\label{corout1}
{\it Let $\gG=(N, M, \alpha, \beta, \Gamma, T, T', \nu)$ be a measured quantum group, and $(b, \ga)$ an outer action of $\gG$ on a von Neumann algebra $A$; 
\newline
(i) the algebra $A$ is a factor if and only if we have :
\[\{n\in N, \alpha(n)\in Z(M)\}=\{n\in N, \beta(n)\in Z(M)\}=\alpha(N)\cap\hat{\beta}(N)=\mathbb{C}\]
Such a measured quantum groupoid is called connected. Then, the scaling operator of $\gG$ is a scalar. 
\newline
(ii) $A\rtimes_\ga\gG$ is a factor if and only if $\widehat{\gG}$ is connected. 
\newline
(iii) both $A$ and $A\rtimes_\ga \gG$ are factors if and only if both $\gG$ and $\widehat{\gG}$ are connected; then, we shall say that $\gG$ is biconnected. We have then :
\[\alpha(N)\cap\beta(N)=\alpha(N)\cap\hat{\beta}(N)=\mathbb{C}\]}
\begin{proof} 
We clearly have, by definition of $\hat{\beta}$, that :
\[\alpha(N)\cap\hat{\beta}(N)\subset\{n\in N, \alpha(n)\in Z(M)\}=\{n\in N, \alpha(n)=\hat{\beta}(n)\}\subset\alpha(N)\cap\hat{\beta}(N)\]
moreover, using the co-inverse $R$, it is clear that $\{n\in N, \alpha(n)\in Z(M)\}=\mathbb{C}$ is equivalent to $\{n\in N, \beta(n)\in Z(M)\}=\mathbb{C}$; then, the first part of (i) is given by \ref{propout}. In that situation, we get immediately that the scaling operator of $\gG$, which belongs to $\alpha(N)\cap Z(M)$, must be a scalar, which finishes the proof of (i). 
\newline
By applying (i) to the action $\tilde{\ga}$ of $\widehat{\gG}^c$ (which is outer by \ref{dualout}), we get (ii), and (i) and (ii) give (iii). \end{proof}

\subsection{Corollary}
\label{corout}
{\it Let $\gG=(N, M, \alpha, \beta, \Gamma, T, T', \nu)$ be a measured quantum groupoid, such that $\alpha(N)\subset Z(M)$, and $(b, \ga)$ an outer action of $\gG$ on a von Neumann algebra $A$; then, we have $Z(A)=b(N)$; let us write $N=L^\infty(X, \nu)$; the von Neumann algebra $A$ is decomposable and can be written $A=\int_X^\oplus A^xd\nu(x)$, and, for $\nu$ almost all $x\in X$, the algebras $A^x$ are factors.}
\begin{proof}
If $\alpha(N)\subset Z(M)$, we have, using the co-inverse $R$, that $\beta(N)\subset Z(M)$, and, therefore, using \ref{propout}, $Z(A)\subset b(N)$, and, as $\ga(b(n))=1\underset{N}{_b\otimes_\alpha}\beta(n)$ commutes with $A\underset{N}{_b*_\alpha}M$, and, therefore, with $\ga(A)$, we get that $b(N)\subset Z(A)$, which finishes the proof. \end{proof}

\subsection{Example}
\label{ex}
Let $M_0\subset M_1$ a depth 2 inclusion, equipped with a regular operator-valued weight $T_1$ from $M_1$ onto $M_0$, and a normal semi-finite faithful weight $\chi$ on $M'_0\cap M_1$, invariant under $\sigma_t^{T_1}$; let us use all notations of \ref{d2}. There exists a measured quantum groupoid $\gG_2$ and an action $\ga$ of $\gG_2$ on $M_1$, with $M_0=M_1^\ga$. Then, this action $\ga$ is outer: 
in fact, the crossed-product $M_1\rtimes_\ga\gG_2$ is isomorphic with $M_2$, and this isomorphism, described in (\cite{EV}7.6), sends $\ga(M_1)$ on $M_1$, and $1\underset{M'_1\cap M_2}{_b\otimes_\alpha}\hat{\alpha}(n)$ on $n$, for any $n\in M'_1\cap M_2$. Here $b$ is the restriction of the mirroring $j_1$ to $M'_1\cap M_2$, which sends the basis $M'_1\cap M_2$ on $M'_0\cap M_1$, $\alpha$ is the injection of $M'_1\cap M_2$ into $M'_1\cap M_3$, and $\hat{\alpha}$ is the restriction to $M'_1\cap M_2$ of the standard representation of $M'_0\cap M_2$. 

\subsection{Example}
\label{lcqg}
(i) Let $\bf{G}$ be a locally compact quantum group; then an action of $\bf{G}$ is outer (in the sense of \ref{defout}) if and only if it is strictly outer in the sense of Vaes (\cite{V2} 2.5). 
\newline
(ii) let $\gG_i$ be a family of measured quantum groupoids, and $(b_i, \ga_i)$ an action of $\gG_i$ on a von Neumann algebra $A_i$. Let us construct $\gG=\oplus_{i\in I}\gG_i$ (\ref{exMQG}(v)); then, let us define $b=\oplus_{i\in I}b_i$, which will be an injective $*$-antihomomorphism from $\oplus_{i\in I}N_i$ into $\oplus_{i\in I}A_i$, and $\ga=\oplus_{i\in I}\ga_i$, which will be an injective $*$-homomorphism from $\oplus_{i\in I}A_i$ into $\oplus_{i\in I}(A_i\underset{N_i}{_{b_i}*_{\alpha_i}}M_i)=(\oplus_{i\in I}A_i)\underset{N}{_b*_\alpha}M$, where $\alpha=\oplus_{i\in I}\alpha_i$ and $M=\oplus_{i\in I}M_i$; then $(b, \ga)$ is an action of $\gG$ on $\oplus_{i\in I}A_i$, and this action is outer if and only if all the actions $\ga_i$ are outer.

\subsection{Theorem}
\label{thgeo}
{\it Let $\gG=(N, M, \alpha, \beta, \Gamma, T, T', \nu)$ be a measured quantum groupoid; then, are equivalent :
\newline
(i) there exists a depth 2 inclusion $M_0\subset M_1$, equipped with a regular operator-valued weight $T_1$ from $M_1$ onto $M_0$, and a normal semi-finite faithful weight $\chi$ on $M'_0\cap M_1$, invariant under $\sigma_t^{T_1}$, such that $\gG=\gG(M_0\subset M_1)$; 
\newline
(ii) there exists a von Neumann algebra $A$, and $(b, \ga)$ a outer action of $\gG$ on $A$. }

\begin{proof}
Let us suppose (i); let $M_0\subset M_1\subset M_2\subset ...$ be Jones' tower associated to the inclusion $M_0\subset M_1$; then (\ref{d2}), $\widehat{\gG}^o=\gG(M_1\subset M_2)$, and, therefore, $\gG^{oc}=\gG(M_2\subset M_3)$, and $\gG=\gG(M_4\subset M_5)$. Applying \ref{ex} to the inclusion $M_3\subset M_4$, we get (ii). 
\newline
Let us suppose (ii); then, using \ref{thout}, we have $\gG^{oc}=\gG(\ga)=\gG(\ga(A)\subset A\rtimes_\ga\gG)$. Using (\cite{E5} 13.9), we get that Jones'tower associated to the incusion $\ga(A)\subset A\rtimes_\ga\gG$ is :
\[\mathfrak a(A)\underset{N^o}{_{\hat{\alpha}}\otimes_\beta}1\subset \tilde{\mathfrak a}(A\rtimes_\mathfrak a\gG)\subset (A\rtimes_\mathfrak a\gG)\rtimes_{\tilde{\mathfrak a}}\widehat{\gG}^c\subset (A\rtimes_\mathfrak a\gG)\underset{N^o}{_{\hat{\alpha}}*_\beta}\mathcal L(H_\Phi)\]
and, therefore, we get that $\gG=\gG[(A\rtimes_\mathfrak a\gG)\rtimes_{\tilde{\mathfrak a}}\widehat{\gG}^c\subset (A\rtimes_\mathfrak a\gG)\underset{N^o}{_{\hat{\alpha}}*_\beta}\mathcal L(H_\Phi)]$, which gives (i). \end{proof}

\section{Faithful actions of a measured quantum groupoid}
\label{faithful}
In this chapter, we define faithful actions of a measured quantum groupoid (\ref{minimal}), and we prove some links between faithful and outer actions (\ref{propminimal}, \ref{thouter}). 
\subsection{Definition}
\label{minimal}
Let $\gG=(N, M, \alpha, \beta, \Gamma, T, T', \nu)$ be a measured quantum groupoid, and $(b, \ga)$ an action of $\gG$ on a von Neumann algebra $A$; we shall say that the action $(b, \ga)$ is faithful if 
\[\{(\omega_\eta\underset{N}{_b*_\alpha}id)\ga(x), \eta\in D(L^2(A)_b, \nu^o), x\in A\}"=M\]
We shall say that $(b, \ga)$ is minimal if it is faithful and if $A\cap (A^{\ga})'=b(N)$. 

\subsection{Example}
\label{exfaithful}
Let $\gG=(N, M, \alpha, \beta, \Gamma, T, T', \nu)$ be a measured quantum groupoid, $(\beta, \Gamma)$ the action of $\gG$ on $M$ defined in (\cite{E5}, 6.10). Using \cite{E5} 3.6 (ii) and 3.8 (vii), we get that the von Neumann algebra generated by the set $\{(\omega_\eta\underset{N}{_\beta*_\alpha}id)\Gamma (x), \eta\in D(H_\beta, \nu^o), x\in M\}$ is equal to $M$, which says that $(\beta, \Gamma)$ is faithful. 

\subsection{Proposition}
\label{propfaithful}
{\it Let $\gG=(N, M, \alpha, \beta, \Gamma, T, T', \nu)$ be a measured quantum groupoid, and $(b, \ga)$ an action of $\gG$ on a von Neumann algebra $A$; let $A_1$ be a von Neumann subalgebra of $A$ such that $b(N)\subset A_1\subset A$, and such that $\ga(A_1)\subset A_1\underset{N}{_b*_\alpha}M$; therefore $(b, a_{|A_1})$ is an action of $\gG$ on $A_1$; moreover, if $(b, a_{|A_1})$ is faithful, then $(b, \ga)$ is faithful. }
\begin{proof}
Trivial. \end{proof}

\subsection{Proposition}
\label{propminimal}
{\it Let $\gG=(N, M, \alpha, \beta, \Gamma, T, T', \nu)$ be a measured quantum groupoid, and $(b, \ga)$ a minimal action of $\gG$ on a von Neumann algebra $A$; then $(b, \ga)$ is an outer action. }
\begin{proof}
Let $z\in A\underset{N}{_b*_\alpha}\mathcal L(H)\cap \ga(A)'$; then, using \ref{minimal} and \ref{spatial} , $z$ belongs to :
\begin{multline*}
A\underset{N}{_b*_\alpha}\mathcal L(H)\cap (A^\ga\underset{N}{_b\otimes_\alpha} 1_H)'=
A\underset{N}{_b*_\alpha}\mathcal L(H)\cap (A^\ga)'\underset{N}{_b*_\alpha}\mathcal L(H)\\=
A\cap (A^\ga)'\underset{N}{_b*_\alpha}\mathcal L(H)
=b(N)\underset{N}{_b*_\alpha}\mathcal L(H)=1\underset{N}{_b\otimes_\alpha}\alpha(N)'
\end{multline*}
So, there exists $y\in\alpha(N)'$ such that $z=1\underset{N}{_b\otimes_\alpha}y$. But, as $z$ commutes with $\ga(A)$, we get that $y$ commutes with all elements of the form $(\omega_\eta\underset{N}{_b*_\alpha}id)\ga(x)$, for all $\eta\in D(L^2(A)_b, \nu^o)$ and $x\in A$. Therefore, by \ref{minimal}, we get that $y\in M'$, which finishes the proof, by \ref{thout}. 
\end{proof}


\subsection{Definition}
\label{trivial}
Let $\gG$ be a measured quantum groupoid, $A$ be a von Neumann algebra (with separable predual), and $\theta$ be a normal faithful state on $A$; let us denote $id_N$ the canonical anti-homomorphism from $N$ into $N^o\otimes B$, and $id_A$ the identity of $A$; then, as the fiber product $(A\otimes N^o)*_\alpha M$ can be identified with $A\otimes (M\cap \alpha(N)')$, we get that $(id_N, id_A\otimes\beta)$ is an action of $\gG$ on $A\otimes N^o$, we shall call the trivial action of $\gG$ on $A\otimes N^o$; this generalizes the example (\cite{E5}, 6.2). If $\theta$ denotes a faithful state on $A$, we shall denote $E_\theta$ the normal faithful conditional expectation from $A\otimes N^o$ onto $N^o$ given by the slice map $\theta\otimes id_N$; this conditional expectation satisfies $(E_\theta\underset{N}{_{id_N}*_\alpha}id)(id_A\otimes\beta)=(id_A\otimes\beta)\circ E_\theta$ and is therefore invariant by the trivial action in the sense of \ref{crossed}. Moreover, we have, trivially $(A\otimes N^o)^{id\otimes\beta}=A\otimes \mathbb{C}$.

\subsection{Proposition}
\label{amalgamated}
{\it Let $\gG$ be a measured quantum groupoid; for $i=(1,2)$, let $(b_i, \ga_i)$ be an action of $\gG$ on a von Neumann algebra $\mathcal A_i$; let us suppose that there exists a normal faithful conditional expectation $E_i$ from $\mathcal i$ onto $b_i(N)$, invariant under $\ga_i$, i.e. (\cite{E7}, 7.6 and 7.7, recalled in \ref{crossed}) such that $(E_i\underset{N}{_{b_i}*_\alpha}id)\ga_i=\ga_i\circ E_i$. Then, there exists an action $(b, \ga)$ of $\gG$ on the amalgamated free product $\mathcal A_1\underset{N^o}{\bigstar}\mathcal A_2$, as defined in (\cite{U},2), where $b$ is the anti-isomorphism from $N$ into the canonical subalgebra of $\mathcal A_1\underset{N^o}{\bigstar}\mathcal A_2$ isomorphic to $N^o$, and $\ga$ is given by the composition of the isomorphism $\ga_1\underset{N^o}{\bigstar}\ga_2$ from $\mathcal A_1\underset{N^o}{\bigstar}\mathcal A_2$ onto  $\ga_1(\mathcal A_1)\underset{N^o}{\bigstar}\ga_2(\mathcal A_2)$ constructed, as (\cite{U}, p.366), using (\cite{U}, 2.5), and the inclusion :
\[\ga_1(\mathcal A_1)\underset{N^o}{\bigstar}\ga_2(\mathcal A_2)\subset (\mathcal A_1\underset{N^o}{\bigstar}\mathcal A_2)\underset{N}{_b*_\alpha}M\]
which is given by the formulae $(E_i\underset{N}{_{b_i}*_\alpha}id)\ga_i=\ga_i\circ E_i$. For any $x_i\in\mathcal A_i$, considered as a subalgebra of $\mathcal A_1\underset{N^o}{\bigstar}\mathcal A_2$, we have $\ga(x_i)=\ga_i(x_i)$. }

\begin{proof}
The construction of the application $\ga$ is an application of (\cite{U}2.5). Then, it is straightforward to get it is an action by verifying it on each $\mathcal A_i$. \end{proof}

\subsection{Definition}
\label{amalgamated2}
Let $\gG$ be a measured quantum groupoid, and let $(b, \ga)$ be a faithful action of $\gG$ on a von Neumann algebra $\mathcal A$; let us suppose that there exists a normal faithful conditional expectation $E$ from $\mathcal A$ onto $b(N)$, invariant under $\ga$. Moreover, let $A$ be a von Neumann algebra with separable predual, and $\theta$ a normal faithful state on $A$, and let us consider the trivial action on $\gG$ on $A\otimes N^o$, as defined in \ref{trivial}. 
\newline
Let us construct now the action $\ga_1$ of $\gG$ on the amalgamated free product $\mathcal A\underset{N^o}{\bigstar}(A\otimes N^o)$ of $\mathcal A$ over its subalgebra $b(N)$ with $(A\otimes N^o)$ over its subalgebra $N^o$, constructed as in \ref{amalgamated} using the normal faithful conditional expectation $E$ from $\mathcal A$ onto $b(N)$ and the normal faithful conditional expectation $E_\theta$ from $A\otimes N^o$ onto $N^o$ defined in \ref{trivial}, which are invariant, respectively, towards the action $\ga$ and the trivial action. 
\newline
As the action $\ga$ is a restriction of $\ga_1$ and is faithful, we get that the action $\ga_1$ is faithful also. Moreover, we have trivially $A\otimes\mathbb{C}\subset [\mathcal A\underset{N^o}{\bigstar}(A\otimes N^o)]^{\ga_1}$.

\subsection{Theorem}
\label{thouter}
{\it Let $\gG$ be a measured quantum groupoid; let us suppose that $\gG$ has a faithful action on a von Neumann algebra $\mathcal A$, such that there exists a normal faithful conditional expectation $E$ from $\mathcal A$ onto $b(N)$, invariant under $\ga$; then $\gG$ has an outer action.}

\begin{proof}
Let's use Barnett's result (\cite{Ba}, th. 2); using the notations of \ref{amalgamated}, let us take $A=(A_1, \theta_1)\bigstar (A_2, \theta_2)$, each $\theta_i$ being a faithful state on the von Neumann algebra $A_i$, and let us suppose that there exists $a$ in the centralizer $A_1^{\theta_1}$ such that $\theta_1(a)=0$, and $b$, $c$ in the centralizer $A_2^{\theta_2}$ such that $\theta_2(b)=\theta_2(c)=\theta_2(b^*c)=0$; let's use the normal faithful conditional expectations $(\theta_1\otimes id)$ from $A_1\otimes N^o$ onto $N^o$ and $(\theta_2\otimes id)$ from $A_2\otimes N^o$ onto $N^o$; it is straightforward to get that the amalgated free product $(A_1\otimes N^o)\underset{N}{\bigstar}(A_2\otimes N^o)$ is equal to $A\otimes N^o$, which, by the associativity of the amalgamated free product, leads to :
\[\mathcal A\underset{N}{\bigstar}(A\otimes N^o)=[\mathcal A\underset{N}{\bigstar}(A_1\otimes N^o)]\underset{N}{\bigstar}(A_2\otimes N^o)\]
Then, we get that the elements $a\otimes 1$, $b\otimes 1$, $c\otimes 1$ satisfie the conditions of (\cite{U}, condition I-A of Appendix I), which leads to (\cite{U}, Prop. I-C of Appendix I) :
\[\mathcal A\underset{N}{\bigstar}(A\otimes N^o)\cap \{a\otimes 1, b\otimes 1, c\otimes 1\}'=N^o\]
from which we get :
\begin{eqnarray*}
\mathcal A\underset{N^o}{\bigstar}(A\otimes N^o)\cap [\mathcal A\underset{N^o}{\bigstar}(A\otimes N^o)]^{\ga_1}
&\subset&
\mathcal A\underset{N^o}{\bigstar}(A\otimes N^o)\cap (A\otimes\mathbb{C})'\\
&\subset&
\mathcal A\underset{N^o}{\bigstar}(A\otimes N^o)\cap \{a\otimes 1, b\otimes 1, c\otimes 1\}'\\
&=&
N^o
\end{eqnarray*}
which proves that the action $\ga_1$ is minimal, in the sense of \ref{minimal}, and, therefore, outer, using \ref{propminimal}. \end{proof}

\section{Any measured quantum groupoid has an outer action}
\label{anyout}

In this chapter, following the strategy of \cite{V2}, we construct a faithful action of a measured quantum groupoid $\gG$ on a von Neumann algebra acting on a "relative Fock space". We prove that this action leaves invariant some conditional expectation; then, using \ref{thouter}, inspired again by \cite{V2}, we construct a strictly outer action on some amagamated free product. 

\subsection{Definition and notations}
\label{defnot1}
For any $n\in \mathbb{N}$, let us write $H^{(n)}$ for $L^2(N)$(that we shall identify with the Hilbert space $H_\nu$ given by the G.N.S. construction made from the weight $\nu$) if $n=0$, for $H$ if $n=1$, for $H\underset{\nu^o}{_\alpha\otimes_{\hat{\beta}}}H$ for $n=2$, and, if $n\geq 3$, for the relative tensor product ($n$-times) $H\underset{\nu^o}{_\alpha\otimes_{\hat{\beta}}}H_\alpha....._{\hat{\beta}}H$.
\newline
Each of these Hilbert spaces is equipped with a surjective involutive antilinear isometry, $J_\nu$ on $H^{(0)}$, $J$ on $H^{(1)}$, $\sigma_{\nu}(J\underset{\nu^o}{_\alpha\otimes_{\hat{\beta}}}J)$ on $H^{(2)}$, and $\Sigma_n(J\underset{\nu^o}{_\alpha\otimes_{\hat{\beta}}}J_\alpha ....._{\hat{\beta}}J)$ on $H^{(n)}$ where $\Sigma_n$ means $\Sigma_n (\xi_1\underset{\nu}{_{\hat{\beta}}\otimes_\alpha}\xi_2{}_{\hat{\beta}}{} .... _\alpha\xi_n)=\xi_n\underset{\nu}{_\alpha\otimes_{\hat{\beta}}}\xi_{n-1}{}_\alpha ..... _{\hat{\beta}} \xi_1$. 
\newline
Let us write $\mathcal F(H)=\oplus_n H^{(n)}$, and let $\mathcal J$ be the surjective involutive antilinear isometry constructed by taking the direct sum of all these isometries on $H^{(n)}$. 
\newline
Let us consider the canonical representation of $N$ on $H^{(0)}$, the representation $\alpha$ on $H^{(1)}$, the representation $1\underset{N^o}{_\alpha\otimes_{\hat{\beta}}}\alpha$ on $H^{(2)}$, and the representations $\alpha_n=1\underset{N^o}{_\alpha\otimes_{\hat{\beta}}}1_\alpha .... _{\hat{\beta}}\alpha$ on $H^{(n)}$, and let us write $a$ for the direct sum of all these representations, which is a normal faithful representation of $N$ on $\mathcal F(H)$. 
\newline
Writing $b(n)=\mathcal J a(n^*)\mathcal J$, we construct a normal faithful antirepresentation of $N$ on $\mathcal F(H)$; we easily get that $b$ is, on $H^{(0)}$, equal to the canonical antirepresentation of $N$, that $b$ is, on $H^{(1)}$, equal to $\hat{\beta}$, that $b$ is, on $H^{(2)}$, equal to $\hat{\beta}\underset{N^o}{_\alpha\otimes_{\hat{\beta}}}1$, and, on $H^{(n)}$, equal to $\hat{\beta}_n=\hat{\beta}\underset{N^o}{_\alpha\otimes_{\hat{\beta}}}1_\alpha .... _{\hat{\beta}} 1$.  
\newline
For any $\xi\in D(_\alpha H, \nu)$, let us define on $\mathcal F(H)$ bounded operator $l(\xi)$ by :
\newline
- for any $n\in N$, $l(\xi)\Lambda_\nu(n)=\alpha(n)\xi$;
\newline
- for any $\eta\in H$, $l(\xi)\eta=\xi\underset{\nu^o}{_\alpha\otimes_{\hat{\beta}}}\eta$;
\newline
- for any $\xi_1\underset{\nu^o}{_\alpha\otimes_{\hat{\beta}}}\xi_2{}_\alpha{} .... _{\hat{\beta}}\xi_n\in H^{(n)}$, $l(\xi)(\xi_1\underset{\nu^o}{_\alpha\otimes_{\hat{\beta}}}\xi_2{}_\alpha{} .... _{\hat{\beta}}\xi_n)=\xi\underset{\nu^o}{_\alpha\otimes_{\hat{\beta}}}\xi_1\underset{\nu^o}{_\alpha\otimes_{\hat{\beta}}}\xi_2{}_\alpha{} .... _{\hat{\beta}}\xi_n$.
\newline
Then, we get that $l(\xi)$ belongs to $a(N)'$, and, for $\xi$, $\xi'$ in $D(_\alpha H, \nu)$, we have $l(\xi')^*l(\xi)=b(\langle\xi, \xi'\rangle_{\alpha, \nu})$. 
\newline
We can easily check that $l(\xi)l(\xi)^*$ is equal to $0$ on $H^{(0)}$, is equal to $\theta^{\alpha, \nu}(\xi, \xi)$ on $H^{(1)}$, and to $\theta^{\alpha, \nu}(\xi, \xi)\underset{N^o}{_\alpha\otimes_{\hat{\beta}}}1$ on $H^{(n)}$. Therefore, if $(\xi_i)_{i\in I}$ is an orthogonal $(\alpha, \nu)$ basis of $H$, in the sense recalled in \ref{spatial}, we get that $\sum_i l(\xi_i)l(\xi_i)^*=1-P_{H^{(0)}}$. 
\newline
Let us write $\mathcal A$ for the von Neumann algebra generated by all the operators $l(\xi)$, for $\xi\in D(_\alpha H, \nu)$. From these remarks, we infer that $b(N)\subset \mathcal A\subset a(N)'$, and that $P_{H^{(0)}}\in\mathcal A$. Taking the final support of $l(\xi)P_{H^{(0)}}$, we get that $\theta^{\alpha, \nu}(\xi, \xi)P_{H^{(1)}}$ belongs to $\mathcal A$, and taking again an $(\alpha, \nu)$-orthogonal basis of $H$, we get that $P_{H^{(1)}}$ belongs to $\mathcal A$. By recurence, we get that, for all $n\in\mathbb{N}$, $P_{H^{(n)}}$ belongs to $\mathcal A$. 
\newline
For any $\eta\in D(H_{\hat{\beta}}, \nu^o)$, we define on $\mathcal F(H)$ a bounded operator $r(\eta)$ by :
\newline
- for any $n\in N$, $r(\eta)\Lambda_\nu(n)=\hat{\beta}(n^*)\eta$;
\newline
- for any $\xi\in H$, $r(\eta)\xi=\xi\underset{\nu^o}{_\alpha\otimes_{\hat{\beta}}}\eta$; 
\newline
- for any $\xi_1\underset{\nu^o}{_\alpha\otimes_{\hat{\beta}}}\xi_2{}_\alpha{} .... _{\hat{\beta}}\xi_n\in H^{(n)}$, $r(\eta) (\xi_1\underset{\omega^o}{_\alpha\otimes_{\hat{\beta}}}\xi_2{}_\alpha{} .... _{\hat{\beta}}\xi_n)=\xi_1\underset{\omega^o}{_\alpha\otimes_{\hat{\beta}}}\xi_2{}_\alpha{} .... _{\hat{\beta}}\xi_n\underset{\omega^o}{_\alpha\otimes_{\hat{\beta}}}\eta$. 
\newline
Then, we easily get that $r(\eta)=\mathcal J l(J\eta)\mathcal J$, and that $r(\eta)\in\mathcal A'$, from which we get that $\mathcal J\mathcal A\mathcal J\subset\mathcal A'$. 
\newline
Let us now consider a faithful normal state $\omega$ on $N$ and the G.N.S. construction ($H_\omega$, $\pi_\omega$, $\Lambda_\omega(1)$) made from $\omega$. There exists a unique unitary $u$ from $H_\omega$ onto $H_\nu=L^2(N)$, such that $u^*nu=\pi_\omega(n)$, for all $n\in N$, and $uJ_\omega=J_\nu u$; then, for any $p\in N$, analytic with respect to $\nu$, using \ref{spatial} and these properties of $u$, we have :
\begin{eqnarray*}
l(\alpha(p)\xi)u\Lambda_\omega(1)
&=&
l(\xi)J_\nu\sigma^\nu_{-i/2}(p^*)J_\nu u\Lambda_\omega(1)\\
&=&
l(\xi)uJ_\omega\pi_\omega(\sigma^\nu_{-i/2}(p^*))J_\omega\Lambda_\omega(1)\\
\end{eqnarray*}
Using the weak density of the analytic elements in $N$, we get that the closure of $l(\mathcal A)u\Lambda_\omega(1)$ contains, for any $\xi\in D(_\alpha H, \nu)$, the subspace $l(\xi)u\pi_\omega(N)'\Lambda_\omega(1)$; therefore, it contains $l(\xi)L^2(N)$, and, by definition, it contains $\xi$. 
\newline
On the other hand, we have $J_\nu u\Lambda_\omega (1)=uJ_\omega\Lambda_\omega(1)=u\Lambda_\omega(1)$; in the sequel, we shall skip the unitary $u$, and consider the vector $\Lambda_\omega(1)$ as an element of $L^2(N)$, invariant by $J_\nu$. 

\subsection{Proposition}
\label{propA}
{\it Let's take the notations of \ref{defnot1}; we have :
\newline
(i) the state $\Omega(X)=(X\Lambda_\omega(1)|\Lambda_\omega(1))$ on $\mathcal A$ is faithful;
\newline
(ii) let $E(X)=b(\langle X\Lambda_\omega (1), \Lambda_\omega(1)\rangle_{a, \omega})$; then, $E$ is a normal faithful conditional expectation from $\mathcal A$ onto $b(N)$. }

\begin{proof}
We had seen in \ref{defnot1} that any $\xi\in D(_\alpha H, \nu)$ belongs to $\overline{\mathcal A \Lambda_\omega(1)}$; therefore,  $\overline{\mathcal A \Lambda_\omega(1)}$ contains $H$; the same way, we get that, for all $n\in\mathbb{N}$,  ${\overline{\mathcal A H^{(n)}}}$ contains $H^{(n+1)}$, and, therefore, the vector $\Lambda_\omega (1)$ is cycling for $\mathcal A$; as $\mathcal J \Lambda_\omega (1)=\Lambda_\omega (1)$, by \ref{defnot1} again, we see that this vector is cycling also for $\mathcal J \mathcal A \mathcal J$, and, therefore, also for $\mathcal A '$; so, $\Lambda_\omega (1)$ is separating for $\mathcal A$, from which we get (i). 
\newline
Let us write, for $X\in\mathcal A$, $E(X)=b(\langle X\Lambda_\omega (1), \Lambda_\omega(1)\rangle_{a, \omega})$; $E$ is a positive bounded application from $\mathcal A$ on $b(N)$; moreover, for any $n\in N$, we have $b(n)\Lambda_\omega(1)=J_\omega n^*J_\omega\Lambda_\omega(1)=\sigma^\omega_{-i/2}(n)\Lambda_\omega(1)$, and $R^{a, \omega}(b(n)\Lambda_\omega(1))=R^{a, \omega}(\Lambda_\omega (1))J_\omega n^*J_\omega$; so, we get that $E(b(n))=b(n)$, and, therefore $E^2=E$, and $E$ is a conditional expectation. As $\Omega(X)=\Omega\circ E(X)$, we get, using (i), that $E$ is faithful, which finishes the proof. \end{proof}

\subsection{Proposition}
\label{propcorep}
{\it Let $\gG=(N, M, \alpha, \beta, \Gamma, T, T', \nu)$ be a measured quantum groupoid; let us use the notations of \ref{defnot1}. Then :
\newline
(i) $\sigma_{\nu^o}W\sigma_{\nu^o}$ is a corepresentation of $\gG$ on the $N-N$ bimodule $_\alpha H_{\hat{\beta}}$; 
\newline
(ii) there exists a unitary 
\[(\sigma_{\nu^o}W\sigma_{\nu^o})_{1,n}(\sigma_{nu^o}W\sigma_{\nu^o})_{2,n}...(\sigma_{\nu^o}W\sigma_{\nu^o})_{n-2,n}(\sigma_{\nu^o}W\sigma_{\nu^o})_{n-1,n}\] from $H^{(n-1)}\underset{\nu}{_{\alpha_n}\otimes_\beta}H$ to $H^{(n-1)}\underset{\nu^o}{_{\hat{\beta}_n}\otimes_\alpha}H$, which is a corepresentation of $\gG$ on the $N-N$ bimodule $_{{\hat{\beta}}_n}H_{\alpha_n}^{(n)}$.
\newline
(iii) by taking the sum of all these, we can define a corepresentation $\mathcal F(\sigma_{\nu^o}W\sigma_{\nu^o})$ of $\gG$ on the $N-N$ bimodule $_b\mathcal F(H)_a$. }
\begin{proof}
Result (i) is nothing but (\cite{E5}, 5.6). It is then easy to get, at least formally, that $(\sigma_{\nu^o}W\sigma_{\nu^o})_{1,3}(\sigma_{\nu^o}W\sigma_{\nu^o})_{2,3}$ is, using (\cite{E5}, 5.1), a corepresentation of $\gG$ on the $N-N$ bimodule $_{\hat{\beta}_2}H^{(2)}_{\alpha_2}$, and we can get by recurrence a proof of (ii). The proof of (iii) is then straightforward. 
\end{proof}

\subsection{Theorem}
\label{thfaithful}
{\it Let $\gG=(N, M, \alpha, \beta, \Gamma, T, T', \nu)$ be a measured quantum groupoid; let us use the notations of \ref{defnot1}, \ref{propA} and \ref{propcorep}. Then the corepresentation $\mathcal F(\sigma_{\nu^o}W\sigma_{\nu^o})$ of $\gG$ on the $N-N$ bimodule $_b\mathcal F(H)_a$ implements, in the sense of (\cite{E5}, 6.6) an action $(b, \ga)$ of $\gG$ on $\mathcal A$ defined, for all $X\in A$ by :
\[\ga(X)=\mathcal F(\sigma_{\nu^o}W\sigma_{\nu^o})(X\underset{N^o}{_a\otimes_\beta}1)\mathcal F(\sigma_{\nu^o}W_\omega\sigma_{\nu^o})^*\]
Moreover, this action is faithful, and the faithful conditional expectation $E$ is invariant by $\ga$. }

\begin{proof}
Using (\cite{E5}, 6.6), we get that $\mathcal F(\sigma_{\nu^o}W\sigma_{\nu^o})$ implements an action on $a(N)'$. Moreover, we get, for $\xi$ and $\eta$ in $D(_\alpha H, \nu)$, and $\eta'\in D(_\alpha H, \nu)\cap D(_\beta H, \nu^o)$, that :
\[(id\underset{N}{_b*_\alpha}\omega_{\eta, \eta'})[\mathcal F(\sigma_{\nu^o}W\sigma_{\nu^o})(l(\xi)\underset{N^o}{_a\otimes_\beta}1)\mathcal F(\sigma_{\nu^o}W\sigma_{\nu^o})^*]=l[(i*\omega_{\eta, \eta'})(\sigma_{\nu^o}W\sigma_{\nu^o})\xi]\]
From which we get that 
\[\mathcal F(\sigma_{\nu^o}W\sigma_{\nu^o})(\mathcal A\underset{N^o}{_a\otimes_\beta}1)\mathcal F(\sigma_{\nu^o}W\sigma_{\nu^o})^*\subset \mathcal A\underset{N}{_b*_\alpha}M\]
which gives that $(b, \ga)$ is an action of $\gG$ on $\mathcal A$. 
\newline
Moreover, using the formula :
\[(id\underset{N}{_b*_\alpha}\omega_{\eta, \eta'})\ga(l(\xi))=l((id*\omega_{\eta, \eta'})(\sigma_{\nu^o}W\sigma_{\nu^o})\xi)\]
we get that, for any $\zeta\in D(H_{\hat{\beta}}, \nu^o)$ and $n\in\gN_\nu$, we get :
\begin{eqnarray*}
((\omega_{\Lambda_\nu(n), \zeta}\underset{N}{_b*_\alpha}id)\ga(l(\xi))\eta|\eta')
&=&
(l((id*\omega_{\eta, \eta'})(\sigma_{\nu^o}W\sigma_{\nu^o})\xi)\Lambda_\nu(n)|\zeta)\\
&=&
(\alpha(n)(id*\omega_{\eta, \eta'})(\sigma_{\nu^o}W\sigma_{\nu^o})\xi|\zeta)\\
&=&
((i*\omega_{\alpha(n)\eta, \eta'})(\sigma_{\nu^o}W\sigma_{\nu^o})\xi|\zeta)\\
&=&
(\omega_{\xi, \zeta}*id)(\sigma_{\nu^o}W\sigma_{\nu^o})\alpha(n)\eta|\eta')
\end{eqnarray*}
from which we get that :
\[(\omega_{\Lambda_\nu(n), \zeta}\underset{N}{_b*_\alpha}id)\ga(l(\xi))=(\omega_{\xi, \zeta}*id)(\sigma_{\nu^o}W\sigma_{\nu^o})\alpha(n)\]
So, using \ref{MQG}, we get that the weak closure of all the elements of the form $(\omega_\eta\underset{N}{_b*_\alpha}id)\ga(x)$, for $\eta\in D(L^2(A)_b, \nu^o)$ and $x\in A$, contains all elements in $M$, and, therefore, this action $\ga$ is faithful. 
Finally, we have, for any $X\in\mathcal A$ :
\[(\Omega\underset{N}{_b*_\alpha}id)\ga(X)=\beta(\langle X\Lambda_\omega(1), \Lambda_\omega(1)\rangle_{\alpha, \omega})\]
and, therefore :
\begin{eqnarray*}
(E\underset{N}{_b*_\alpha}id)\ga(X)
&=&
1\underset{N}{_b\otimes_\alpha}(\Omega\underset{N}{_b*_\alpha}id)\ga(X)\\
&=&
1\underset{N}{_b\otimes_\alpha}\beta(\langle X\Lambda_\omega(1), \Lambda_\omega(1)\rangle_{\alpha, \omega})\\
&=&
\ga(b(\langle X\Lambda_\omega(1), \Lambda_\omega(1)\rangle_{\alpha, \omega}))\\
&=&
\ga(E(X))
\end{eqnarray*}
which finishes the proof. \end{proof}

\subsection{Theorem}
\label{thgeo2}
{\it Let $\gG$ be a measured quantum groupoid; then :
\newline
(i) there exists an outer action of $\gG$; 
\newline
(ii) there exists a depth 2 inclusion $M_0\subset M_1$, equipped with a regular operator-valued weight $T_1$ from $M_1$ onto $M_0$, and a normal semi-finite faithful weight $\chi$ on $M'_0\cap M_1$, invariant under $\sigma_t^{T_1}$, such that $\gG=\gG(M_0\subset M_1)$. }

\begin{proof}
Result (i) is clear by \ref{thouter} and (ii) by \ref{thgeo}. \end{proof}

\section{Outer actions on semi-finite and finite von Neumann algebras}
\label{S}
In this chapter, we study the case when a measured quantum groupoid is acting outerly on a semi-finite von Neumann algebra (\ref{semi-finite}, \ref{rhoinv}, \ref{thtensor}), or a finite von Neumann algebra (\ref{finite}, \ref{finite2}). S. Vaes had proved (\cite{V2}, 3.5) that, if a locally compact quantum group acts outerly on a $II_1$ factor, then its scaling group $\tau_t$ is trivial. Here the situation is much more complicated, as it is known,  since M.-C. David's result (\cite{D}), that any connected finite dimensional measured quantum groupoid (with an antipode which is involutive on the two copies of the basis) acts outerly on the hyperfinite $II_1$ factor (and, instead of the locally compact quantum case), there are finite dimensional quantum groupoids with a non trivial scaling group). 

\subsection{Definition}
\label{defweighted}
Let $\gG$ be a measured quantum groupoid, and let $(b, \ga)$ be an action on a von Neumann algebra $A$. We shall say (\cite{E7}, 4.1) that this action is weighted if there exists a normal, semi-finite faithful operator-valued weight $\gT$ from $A$ onto $b(N)$. Then, the weight $\psi=\nu^o\circ b^{-1}\circ\gT$ will be called lifted from $\nu^o$ (or lifted). Then, for any lifted weight $\psi$ on $A$, it is possible to define a 2-cocycle $(D\psi\circ\ga:D\psi)_t=\Delta_{\tilde{\psi}}^{it}(\Delta_{\psi}^{-it}\underset{N}{_b\otimes_\alpha}\Delta_{\widehat{\Phi}}^{it})$ in $A\underset{N}{_b*_\alpha}(M\cap\beta(N)')$ which satisfies, for all $s$, $t$ in $\mathbb{R}$ (\cite{E7}, 7.2 and 7.3) :
\[(D\psi\circ\ga:D\psi)_{s+t}=(D\psi\circ\ga:D\psi)_s(\sigma_s^\psi\underset{N}{_b*_\alpha}\tau_s)((D\psi\circ\ga:D\psi)_t)\]
\[(id\underset{N}{_b*_\alpha}\Gamma)((D\psi\circ\ga:D\psi)_t)=(\ga\underset{N}{_b*_\alpha}id)((D\psi\circ\ga:D\psi)_t)((D\psi\circ\ga:D\psi)_t\underset{N}{_\beta\otimes_\alpha}1)\]
This 2-cocycle is, by definition, Connes' cocycle derivative $(D\overline{\psi}:D\underline{\psi})_t$, where $\overline{\psi}$ is the bidual weight of $\psi$, defined on $A\underset{N}{_b*_\alpha}\mathcal L(H)$ which is canonically isomorphic to the bicrossed product (\cite{E5}, 11.6), and the weight $\underline{\psi}$ is equal to $\overline{\nu^o}\circ(\psi\underset{N}{_b*_\alpha}id)$, where $\overline{\nu^o}$ is a normal semi-finite faithful weight on $\alpha(N)'$, such that $\frac{d\overline{\nu^o}}{d\nu^o}=\Delta_{\widehat{\Phi}}^{-1}$ (\cite{E7}, 4.6), and $\psi\underset{N}{_b*_\alpha}id$ is a normal semi-finite faithful operator-valued weight from $A\underset{N}{_b*_\alpha}\mathcal L(H)$ onto $\alpha(N)'$ (\cite{E7}, 4.4); moreover, we have then 
$\frac{d\overline{\psi}}{d\psi^o}=\Delta_{\tilde{\psi}}$, where $\tilde{\psi}$ is the dual weight of $\psi$, defined on the crossed-product $A\rtimes_\ga\gG$, and 
$\frac{d\underline{\psi}}{d\psi^o}=\Delta_\psi\underset{N}{_b\otimes_\alpha}\Delta_{\widehat{\Phi}}^{-1}$, which leads to the result.

\subsection{Proposition}
\label{T}
{\it Let $\gG$ be a measured quantum groupoid, and let $(b, \ga)$ be a weighted strictly outer action of $\gG$ on a von Neumann algebra $A$. Let $t\in \mathbb{R}$ be in Connes'invariant $T(A)$ (\cite{St}, 27.1); then, there exists $v\in M\cap\alpha(N)'\cap \beta(N)'$ such that $\Gamma(v)=v\underset{N}{_b\otimes_\alpha}v$, and $\tau_t(x)=vxv^*$, for all $x\in M$; moreover, we have $\sigma_t^\nu=id$.}

\begin{proof}
Let $\psi$ be a lifted weight; as $\sigma_t^\psi$ is interior, there exists a unitary $w\in A$ such that $\sigma_t^\psi(x)=wxw^*$ for all $x\in A$; therefore, we get that $\Delta_\psi^{it}=wJ_\psi wJ_\psi$, and, using (\cite{E7}, 7.1 and 7.2), we get that :
\[(D\psi\circ\ga:D\psi)_t=\Delta_{\tilde{\psi}}^{it}(w^*J_\psi w^*J_\psi\underset{N}{_b\otimes_\alpha}\Delta_{\widehat{\Phi}}^{it})\]
One should note that it is possible to define the unitary $w^*\underset{\nu}{_b\otimes_\alpha}\Delta_{\widehat{\Phi}}^{it}$ on elementary tensors (and then extends it to the Hilbert space $H_\psi\underset{\nu}{_b\otimes_\alpha}H$) because we have, for all $n\in N$, $w^*b(n)w=\sigma_{-t}^\psi(b(n))=b(\sigma_t^\nu(n))$ and $\Delta_{\widehat{\Phi}}^{it}\alpha(n)\Delta_{\widehat{\Phi}}^{-it}=\tau_{t}(\alpha(n))=\alpha(\sigma_{t}^\nu(n))$. 
Moreover, we have :
\begin{eqnarray*}
(D\psi\circ\ga:D\psi)_t(wJ_\psi wJ_\psi\underset{N}{_b\otimes_\alpha}\Delta_{\widehat{\Phi}}^{-it})\ga(x)
&=&
\Delta_{\tilde{\psi}}^{it}\ga(x)\\
&=&
\sigma_t^{\hat{\psi}}(\ga(x))\Delta_{\tilde{\psi}}^{it}\\
&=&
\ga(\sigma_t^\psi(x))\Delta_{\tilde{\psi}}^{it}\\
&=&
\ga(wxw^*)\Delta_{\tilde{\psi}}^{it}
\end{eqnarray*}
which is equal to $\ga(w)\ga(x)\ga(w^*)(D\psi\circ\ga:D\psi)_t(wJ_\psi wJ_\psi\underset{N}{_b\otimes_\alpha}\Delta_{\widehat{\Phi}}^{-it})$.
\newline
From which we get that $\ga(w^*)(D\psi\circ\ga:D\psi)_t(w\underset{N}{_b\otimes_\alpha}\Delta_{\widehat{\Phi}}^{-it})$ (which belongs to $A\underset{N}{_b*_\alpha}\mathcal L(H)$) commutes with $\ga(x)$, for all $x\in M$. Using then \ref{defout} and \ref{thout}, we get that there exists $u\in M'$ such that :
\[\ga(w^*)(D\psi\circ\ga:D\psi)_t(w\underset{N}{_b\otimes_\alpha}\Delta_{\widehat{\Phi}}^{-it})=1\underset{N}{_b\otimes_\alpha}u\]
or :\[(D\psi\circ\ga:D\psi)_t=\ga(w)(w^*\underset{N}{_b\otimes_\alpha}u\Delta_{\widehat{\Phi}}^{it})\]
from which we deduce that $u\Delta_{\widehat{\Phi}}^{it}=v$ belongs to $M\cap\beta(N)'$; so $v\Delta_{\widehat{\Phi}}^{-it}$ belongs to $M'$, and $vxv^*=\Delta_{\widehat{\Phi}}^{it}x\Delta_{\widehat{\Phi}}^{-it}=\tau_{t}(x)$. So, the automorphism $\tau_t$ is interior; as $v\in\beta(N)'$, we get that $\beta(\sigma^\nu_t(n))=\tau_t(\beta(n))=\beta(n)$, which implies that $\sigma_t^\nu=id$. So, we get that $wb(n)w^*=\sigma^\psi_t(b(n))=b(\sigma_{-t}^\nu(n))=b(n)$, and, therefore, $w\in A\cap b(N)'$.
Moreover, as :
\[(D\psi\circ\ga:D\psi)_t=\ga(w)(w^*\underset{N}{_b\otimes_\alpha}v)\]
and $w$ commutes with $\beta(N)$, we get that $v\in\alpha(N)'$; moreover, the cocycle property with respect to $\ga$ gives that :
\begin{eqnarray*}
(id\underset{N}{_b*_\alpha}\Gamma)\ga(w)(w^*\underset{N}{_b\otimes_\alpha}\Gamma(v))
&=&
(\ga\underset{N}{_b*_\alpha}id)\ga(w)(\ga(w^*)\underset{N}{_b\otimes_\alpha}v)(\ga(w)(w^*\underset{N}{_b\otimes_\alpha}v)\underset{N}{_\beta\otimes_\alpha}1)\\
&=&
((\ga\underset{N}{_b*_\alpha}id)\ga(w)(w^*\underset{N}{_b\otimes_\alpha}v\underset{N}{_\beta\otimes_\alpha}v)
\end{eqnarray*}
from which we deduce that $\Gamma(v)=v\underset{N}{_\beta\otimes_\alpha}v$.
\end{proof}

\subsection{Proposition}
\label{semi-finite}
{\it Let $\gG$ be a measured quantum groupoid, and let $(b, \ga)$ be a weighted outer action of $\gG$ on a semi-finite von Neumann algebra $A$; then :
\newline
(i) there exists a positive non singular operator $\rho$ affiliated to $M\cap\alpha(N)'\cap\beta(N)'$, such that, for any $t\in \mathbb{R}$, $x\in M$, we have :
\[\Gamma(\rho)=\rho\underset{N}{_\beta\otimes_\alpha}\rho\]
\[\tau_t(x)=\rho^{it}x\rho^{-it}\]
from which we deduce that $\rho$ commutes with $\Delta_{\widehat{\Phi}}$. 
\newline
(ii) the weight $\nu$ is a trace, and, for any normal semi-finite faithful trace $\theta$ on $A$, $\theta$ is lifted, and $(D\overline{\theta}:D\underline{\theta})_t=1\underset{N}{_b\otimes_\alpha}\rho^{it}$, with the notations of \ref{defweighted}. 
\newline
(iii) we have $Z(A)\subset A^\ga$, $Z(A)\underset{N}{_b\otimes_\alpha}\mathbb{C}\subset Z(A\rtimes_\ga\gG)$, and $\alpha(N)\cap Z(M)\subset \alpha(N)\cap Z(\widehat{M})$. }

\begin{proof}
Let us apply \ref{T} to the hypothesis; we get that $\nu$ is a trace, and, therefore, that any normal semi-finite faithful trace $\theta$ on $A$ is lifted from $\nu^o$; we obtain then, for any $t\in\mathbb{R}$, the existence of a unitary $v_t$ in $M\cap\alpha(N)'\cap\beta(N)'$ such that $\Gamma(v_t)=v_t\underset{N}{_\beta\otimes_\alpha}v_t$, $\tau_t(x)=v_txv_t^*$, for all $x\in M$, and $(D\theta\circ\ga :D\theta)_t=1\underset{N}{_b\otimes_\alpha}v_t$; it is therefore clear that the application $t\mapsto v_t$ is continuous; moreover, the cocycle relation (with respect to $(id\underset{N}{_b*_\alpha}\tau_t)$) of $(D\theta\circ\ga :D\theta)_t$ leads to 
$v_{s+t}=v_s\tau_s(v_t)$.  
\newline
But we get also that, for all $t\in\mathbb{R}$, we have $\Delta_{\tilde{\theta}}^{it}=1\underset{N}{_b\otimes_\alpha}v_t\Delta_{\widehat{\Phi}}^{-it}$. From which we deduce that $t\mapsto v_t\Delta_{\widehat{\Phi}}^{-it}$ is a one-parameter group of unitaries. So, for any $s$, $t$ in $\mathbb{R}$, we have $v_s\Delta_{\widehat{\Phi}}^{-is}v_t\Delta_{\widehat{\Phi}}^{-it}=v_{s+t}\Delta_{\widehat{\Phi}}^{-i(s+t)}$. From the cocycle relation of $v_t$, we then get that $\Delta_{\widehat{\Phi}}^{-is}v_t=\tau_s(v_t)\Delta_{\widehat{\Phi}}^{-is}$, and, therefore, that $\tau_s(v_t)=\tau_{-s}(v_t)$; from which we deduce that the unitaries $v_t$ are invariant under $\tau_s$, and, therefore, that $t\mapsto v_t$ is a one-parameter group of unitaries in $M\cap\alpha(N)'\cap\beta(N)'$. 
From which, with the help of \ref{T}, we finish the proof of (i) and (ii). 
\newline
Moreover, let now $k\in Z(N)^+$, such that $\alpha(k)$ belongs to $Z(M)$; then $\beta(k)=R(\alpha(k))$ belongs also to $Z(M)$, using \ref{propout}, we get that $b(k)\in Z(A)$; Let us write $k=\int_0^{\|k\|}\lambda de_\lambda$, and $k_n=\int_{1/n}^{\| k\|}\lambda de_\lambda$; then $k_n$ is invertible, and, for any normal semi-finite faithful trace $\theta$ on $A$, there exists a normal semi-finite faithful trace $\theta_n$ on $A$ such that $(D\theta_n: D\theta)_t=b(k_n)^{it}$. We then obtain that :
\[(D\overline{\theta_n}:D\overline{\theta})_t=\ga(b(k_n)^{it})=1\underset{N}{_b\otimes_\alpha}\beta(k_n)^{it}\]
and, on the other hand :
\[(D\underline{\theta_n}:D\underline{\theta})_t=b(k_n)^{it}\underset{N}{_b\otimes_\alpha}1=1\underset{N}{_b\otimes_\alpha}\alpha(k_n)^{it}\]
Applying then (ii) to the traces $\theta$ and $\theta_n$, as $\rho$ commutes with $\alpha(N)$ and $\beta(N)$, we get that $\alpha(k_n)=\beta(k_n)$, and, when $n$ goes to $\infty$, $\alpha(k)=\beta(k)$, from which we get that $\alpha(k)$ belongs to $Z(\widehat{M})$. 
\newline
Moreover, let now $x\in Z(A)$; using again \ref{propout}, we get that there exists $k\in Z(N)$ such that $\alpha(k)$ belongs to $Z(M)$, and $x=b(k)$; but, now, we have, as we proved that $\alpha(k)=\beta(k)$ :
\[\ga(x)=1\underset{N}{_b\otimes_\alpha}\beta(k)=1\underset{N}{_b\otimes_\alpha}\alpha(k)=b(k)\underset{N}{_b\otimes_\alpha}1\]
which proves that $Z(A)\subset A^\ga$. On the other hand, as $\beta(k)$ belongs to $Z(M)$, we have also :
\[x\underset{N}{_b\otimes_\alpha}1=\ga(x)=1\underset{N}{_b\otimes_\alpha}\beta(k)=1\underset{N}{_b\otimes_\alpha}\hat{\alpha}(n)\]
and, using again \ref{propout}, we get that $x\underset{N}{_b\otimes_\alpha}1$ belongs to $Z(A\rtimes_\ga\gG)$, which finishes the proof. 

\end{proof}

\subsection{Corollary}
\label{corfactor}
{\it Let $\gG$ be a measured quantum groupoid, and let $(b, \ga)$ be a weighted outer action of $\gG$ on a semi-finite von Neumann algebra $A$; then, if $A\rtimes_\ga\gG$ is a factor, then $A$ is a factor; equivalently, if $\widehat{\gG}$ is connected, then $\gG$ is connected also. }
\begin{proof}
This is clear, using \ref{semi-finite}(iii). \end{proof}

\subsection{Theorem}
\label{rhoinv}
{\it Let $\gG$ be a measured quantum groupoid, and let $(b, \ga)$ be a weighted outer action of $\gG$ on a semi-finite von Neumann algebra $A$; let $\theta$ be a normal semi-finite faithful trace on $A$; then, $\nu$ is a trace, there exists a normal semi-finite faithful operator-valued weight $\mathcal T$ from $A$ onto $b(N)$ such that $\theta=\nu\circ b^{-1}\circ\mathcal T$, and we have, for all $x\in\gN_{\theta}\cap\gN_{\mathcal T}$ :
\[(\theta\underset{N}{_b*_\alpha}id)\ga(x^*x)=\beta\circ b^{-1}\circ\mathcal T(x^*x)\rho^{-1}\]
where $\rho$ is a non singular positive operator affiliated to $M\cap\alpha(N)'\cap\beta(N)'$ had been defined in \ref{semi-finite} and satisfies, for any $t\in \mathbb{R}$, $x\in M$ :}
\[\Gamma(\rho)=\rho\underset{N}{_\beta\otimes_\alpha}\rho\]
\[\tau_t(x)=\rho^{it}x\rho^{-it}\]

\begin{proof}
We had got in \ref{semi-finite} that $\nu$ is a trace, and the existence of the operator $\rho$; as $\rho$ is affiliated to $M\cap\alpha(N)'\cap\beta(N)'$, if we write $\rho=\int_o^\infty \lambda de_\lambda$ and, for all $n\in\mathbb{N}$, $f_n=\int_{1/n}^n de_\lambda$, we get, for any $\xi\in D(_\alpha H, \nu)\cap D(H_\beta, \nu^o)$, that $f_n\xi$ belongs to $D(_\alpha H, \nu)\cap D(H_\beta, \nu^o)\cap\mathcal D(\rho^{-1/2})$, and that $\rho^{-1/2}f_n\xi$ belongs to $D(_\alpha H, \nu)\cap D(H_\beta, \nu^o)$.
\newline
As $(D\overline{\theta}:D\underline{\theta})_t=1\underset{N}{_b\otimes_\alpha}\rho^{it}$, by \ref{semi-finite}, we have, for all $\zeta$, $\zeta'$ in $D(_\alpha H, \nu)$ and $t\in\mathbb{R}$ :
\begin{eqnarray*}
\sigma_t^{\overline{\theta}}(1\underset{N}{_b\otimes_\alpha}\theta^{\alpha, \nu}(\zeta, \zeta'))
&=&
(1\underset{N}{_b\otimes_\alpha}\rho^{it})\sigma_t^{\underline{\theta}}(1\underset{N}{_b\otimes_\alpha}\theta^{\alpha, \nu}(\zeta, \zeta'))(1\underset{N}{_b\otimes_\alpha}\rho^{-it})\\
&=&
1\underset{N}{_b\otimes_\alpha}\rho^{it}\sigma_t^{\overline{\nu^o}}(\theta^{\alpha, \nu}(\zeta, \zeta'))\rho^{-it}\\
&=&
1\underset{N}{_b\otimes_\alpha}\rho^{it}\Delta_{\widehat{\Phi}}^{-it}\theta^{\alpha, \nu}(\zeta, \zeta')\Delta_{\widehat{\Phi}}^{it}\rho^{-it}\\
&=&
1\underset{N}{_b\otimes_\alpha}\theta^{\alpha, \nu}(\rho^{it}\Delta_{\widehat{\Phi}}^{-it}\zeta, \rho^{it}\Delta_{\widehat{\Phi}}^{-it}\zeta')
\end{eqnarray*}
\newline
Let $x\in\gN_\theta$, $\xi\in D(_\alpha H, \nu)\cap D(H_\beta, \nu^o)\cap \mathcal D(\rho^{-1/2})$, $\eta$
in $\mathcal D(\Delta_{\widehat{\Phi}}^{-1/2})\cap \mathcal D(\rho^{-1/2})$, such that $\Delta_{\widehat{\Phi}}^{-1/2}\eta$ and $\rho^{-1/2}\eta$ belong to $D(_\alpha H, \nu)$; we have then :
\[\|\Lambda_\theta(x)\underset{\nu^o}{_a\otimes_\beta}\rho^{-1/2}\xi\underset{\nu^o}{_\alpha\otimes_\beta}J_{\widehat{\Phi}}\Delta_{\widehat{\Phi}}^{-1/2}\eta\|^2
=
\|\Lambda_\theta(x^*)\underset{\nu}{_b\otimes_\alpha}J_{\widehat{\Phi}}\rho^{-1/2}\xi\underset{\nu}{_\beta\otimes_\alpha}\Delta_{\widehat{\Phi}}^{-1/2}\eta\|^2\]
which, using (\cite{E7}, 4.11), is equal to :
\[\|\Lambda_{\overline{\theta}}((1\underset{N}{_b\otimes_\alpha}\theta^{\alpha, \nu}(\Delta_{\widehat{\Phi}}^{-1/2}\eta, \Delta_{\widehat{\Phi}}^{1/2}\rho^{-1/2}\xi))\ga(x^*))\|^2\]
The hypothesis about $\xi$ and $\eta$ give that $(1\underset{N}{_b\otimes_\alpha}\theta^{\alpha, \nu}(\Delta_{\widehat{\Phi}}^{-1/2}\eta, \Delta_{\widehat{\Phi}}^{1/2}\rho^{-1/2}\xi))^*$ belongs to $\mathcal D(\sigma_{-i/2}^{\overline{\theta}})$, and, therfore, we get that :
\[\|\Lambda_\theta(x)\underset{\nu^o}{_a\otimes_\beta}\rho^{-1/2}\xi\underset{\nu^o}{_\alpha\otimes_\beta}J_{\widehat{\Phi}}\Delta_{\widehat{\Phi}}^{-1/2}\eta\|^2
=
\|\Lambda_{\overline{\theta}}(\ga(x)(1\underset{N}{_b\otimes_\alpha}\theta^{\alpha, \nu}(\xi, \rho^{-1/2}\eta))\|^2\]
which, using again the Radon-Nykodym derivative betwwen $\overline{\theta}$ and $\underline{\theta}$, is equal to :
\[
\|\Lambda_{\underline{\theta}}(\ga(x)(1\underset{N}{_b\otimes_\alpha}\theta^{\alpha, \nu}(\xi, \eta))\|^2
=
\overline{\nu^o}(\theta^{\alpha, \nu}(\xi, \eta)^*(\theta\underset{N}{_b*_\alpha}id)\ga(x^*x)\theta^{\alpha, \nu}(\xi, \eta))\]
We had got in (\cite{E7}, 4.11) that $\Lambda_{\overline{\nu^o}}(\theta^{\alpha, \nu}(\xi, \eta))=\xi\underset{\nu^o}{_\alpha\otimes_\beta}J_{\widehat{\Phi}}\Delta_{\widehat{\Phi}}^{-1/2}\eta$, and, therefore, we get :
\[\|\Lambda_\theta(x)\underset{\nu^o}{_a\otimes_\beta}\rho^{-1/2}\xi\underset{\nu^o}{_\alpha\otimes_\beta}J_{\widehat{\Phi}}\Delta_{\widehat{\Phi}}^{-1/2}\eta\|^2
=
((\theta\underset{N}{_b*_\alpha}id)\ga(x^*x)\xi\underset{\nu^o}{_\alpha\otimes_\beta}J_{\widehat{\Phi}}\Delta_{\widehat{\Phi}}^{-1/2}\eta|\xi\underset{\nu^o}{_\alpha\otimes_\beta}J_{\widehat{\Phi}}\Delta_{\widehat{\Phi}}^{-1/2}\eta)\]
from which we infer that :
\begin{multline*}
(\Lambda_{\theta}(x)\underset{\nu^o}{_a\otimes_\beta}\alpha(\langle J_{\widehat{\Phi}}\Delta_{\widehat{\Phi}}^{-1/2}\eta, J_{\widehat{\Phi}}\Delta_{\widehat{\Phi}}^{-1/2}\eta\rangle_{\beta, \nu^o})\rho^{-1/2}\xi|\Lambda_{\theta}(x)\underset{\nu^o}{_a\otimes_\beta}\rho^{-1/2}\xi)=\\
((\theta\underset{N}{_b*_\alpha}id)\ga(x^*x)\alpha(\langle J_{\widehat{\Phi}}\Delta_{\widehat{\Phi}}^{-1/2}\eta, J_{\widehat{\Phi}}\Delta_{\widehat{\Phi}}^{-1/2}\eta\rangle_{\beta, \nu^o})\xi|\xi)
\end{multline*}
which, by density gives, for all $n\in N$ :
\[(\Lambda_{\theta}(x)\underset{\nu^o}{_a\otimes_\beta}\alpha(n)\rho^{-1/2}\xi|\Lambda_{\theta}(x)\underset{\nu^o}{_a\otimes_\beta}\rho^{-1/2}\xi)=((\theta\underset{N}{_b*_\alpha}id)\ga(x^*x)\alpha(n)\xi|\xi)\]
and, therefore $\|\Lambda_{\theta}(x)\underset{\nu^o}{_a\otimes_\beta}\rho^{-1/2}\xi\|^2=((\theta\underset{N}{_b*_\alpha}id)\ga(x^*x)\xi|\xi)$. 
If now $x$ belongs to $\gN_\theta\cap\gN_{\mathcal T}$, $\Lambda_\theta(x)$ belongs to $D(_aH_\theta, \nu)$, and $\langle\Lambda_\theta(x), \Lambda_\theta(x)\rangle_{a, \nu}=b^{-1}(\mathcal T(x^*x))$. We get then that :
\[(\beta\circ b^{-1}(\mathcal T(x^*x))\rho^{-1}\xi|\xi)=((\theta\underset{N}{_b*_\alpha}id)\ga(x^*x)\xi|\xi)\]
and, as we deal on both sides with positive closed operators, by density, we get the result. 
\end{proof}

\subsection{Notations}
\label{notgroupoid}
On the constructive point of view, we shall now prove that any locally compact groupoid acts outerly on a semi-finite von Neumann algebra (\ref{thtensor}); in the case of finite groupoids, this result had been obtained by J.-M. Vallin in (\cite{Val5}, 3.3.11); let's fix the notations : 
let $\mathcal G$ be a locally compact groupoid, in the sense of \cite{R}, equipped with a left Haar system $(\lambda^u)_{u\in\mathcal G^{(0)}}$ and a quasi-invariant measure $\nu$ on the set of units $\mathcal G^{(0)}$. Let us denote $\mu=\int_{\mathcal G^{(0)}}\lambda^x d\nu(x)$. Let us consider the left regular representation $\lambda(g)$ of $\mathcal G$; for any $g\in\mathcal G$, $\lambda(g)$ is a unitary from $L^2(\mathcal G^{r(g)}, \lambda^{r(g)})$ onto $L^2(\mathcal G^{s(g)}, \lambda^{s(g)})$, where, as usual, for any $x\in\mathcal G^{(0)}$, $\mathcal G^x=r^{-1}(x)$. We can consider as well $\lambda^g$ as an orthogonal operator from the real Hilbert space $L_{\mathbb{R}}^2(\mathcal G^{r(g)}, \lambda^{r(g)})$ onto the real Hilbert space $L_{\mathbb{R}}^2(\mathcal G^{s(g)}, \lambda^{s(g)})$; this operator extends to an isomorphism from the Clifford algebra $Cl(L_{\mathbb{R}}^2(\mathcal G^{r(g)}, \lambda^{r(g)}))$ (see \cite{Bla} for details) onto the Clifford algebra $Cl(L_{\mathbb{R}}^2(\mathcal G^{(s(g)}, \lambda^{s(g)}))$;
on each of these algebras, there exists a finite trace, and each GNS representation of these algebras generates a factor, which is the hyperfinite $II_1$ factor if $\mathcal G^{r(g)}$ is infinite (and a finite dimensional factor if it is finite). This construction is just a generalization of \cite{Bla} up to groupoids. 
\newline
So, for any $x\in\mathcal G^{(0)}$, we get a copy $A^x$ of the hyperfinite $II_1$ factor $\mathcal R$ (or a finite dimensional factor), and, for any $g\in\mathcal G$, an isomorphism $\ga_g$ from $A^{r(g)}$ onto $A^{s(g)}$; we obtain then an action $\ga$ of the groupoid $\mathcal G$ on the von Neumann algebra $A=\int_{\mathcal G^{(0)}}^{\oplus}A^xd\nu(x)$ (which is hyperfinite $II_1$); let us write $\tau^x$ for the canonical finite trace on $A^x$. By the unicity of the normalized trace $\tau^x$ on $A^x$, we have, for any $g\in\mathcal G$, and $y$ positive in $A^{r(g)}$, $\tau^{s(g)}(\ga_g(y))=\tau^{r(g)}(x)$. If we consider now $\gG(\mathcal G)$ the canonical abelian measured quantum groupoid associated to $\mathcal G$, we get that the application $\ga$ given by the formula $\ga[\int^{\oplus}_{\mathcal G^{(0)}}a^x d\nu (x)]=\int^{\oplus}_\mathcal G a_g(a^{s(g)})d\mu(g)$ is an action of $\gG(\mathcal G)$ on $A$, together with the isomorphism $b$ of $L^\infty(\mathcal G^{(0)}, \nu)$ with $Z(A)$ (\cite{E5} 6.3); it is clear that the formula $E(\int_{\mathcal G^{(0)}}^{\oplus}a^xd\nu(x))=b(x\mapsto \tau^x(a^x))$ is a normal faithful conditional expectation from $A$ onto $Z(A)$; moreover, we have then $(E\underset{L^\infty(\mathcal G^{(0)}, \nu)}{_b*_\alpha}id)\ga=\ga\circ E$. 

\subsection{Proposition}
\label{proptensor}
{\it Let $\mathcal G$ be a locally compact separable groupoid, and let $\ga$ be a faithful action of $\mathcal G$ on a von Neumann algebra $A=\int^{\oplus}_{\mathcal G^{(0)}}A^x d\nu(x)$, where the algebras $A^x$ are factors, equipped with a faithful state $\omega^x$, invariant by $\ga$, i.e. such that $\omega^{s(g)}\circ\ga_g=\omega^{r(g)}$, for all $g\in\mathcal G$. Let us consider the infinite tensor product $(B^x, \omega_\infty^x)=\otimes_{\mathbb{N}}(A^x, \omega^x)$ of copies of $(A^x, \omega^x)$. It is clear that $(B^x)_{x\in\mathcal G^{(0)}}$ is a continuous field of factors. Then the action   $\tilde{\ga}_g=\otimes_{\mathbb{N}}(\ga_g, \omega^{r(g)})$ defines an outer action of $\mathcal G$ on $B=\int_{\mathcal G^{(0)}}^{\oplus}B^x d\nu(x)$. }

\begin{proof} The proof is completely taken from (\cite{V2}, 5.1), and we shall give the arguments only when it differs. 
\newline
Let's take $a\in B\rtimes_{\tilde{\ga}}\mathcal G\cap \tilde{\ga}(B)'$; as in (\cite{V2}, 5.1), we can prove that $a$ commutes with all elements of the form $1\underset{L^\infty(\mathcal G^{(0)})}{_b\otimes_r}(\omega_\xi\underset{L^\infty(\mathcal G^{(0)})}{_b*_r}id)\ga(x)$, for all $x\in B$, and $\xi\in D(\otimes_{\mathbb{N}}(L^2(A^x), \omega^x)_b, \nu)$, where $b$ means the isomorphism from $L^\infty(\mathcal G^{(0)}, \nu)$ onto $Z(B)$, and $r$ the injection of $L^\infty(\mathcal G^{(0)}, \nu)$ into $L^\infty(\mathcal G, \mu)$ given the the range function. As $\ga$ is faithful, these elements are functions on $\mathcal G$ which separate the points of $\mathcal G$, and we get that the commutant of these elements is equal to $\mathcal L(L^2(B))\underset{L^\infty(\mathcal G^{(0)}, \nu)}{_b*_r}L^\infty(\mathcal G, \mu)$. Using then (\cite{E5}, 9.4 and 11.5), we get that $a$ belongs to $\tilde{\ga}(B)$, and, therefore, to $\tilde{\ga}(Z(B))=1\underset{L^\infty(\mathcal G^{(0)}, \nu)}{_b\otimes_r}s(L^\infty(\mathcal G^{(0)}, \nu)$, where $s$ is the injection of $L^\infty(\mathcal G^{(0)}, \nu)$ into $L^\infty(\mathcal G, \mu)$ given the the source function, which is here equal to $\hat{r}$ (because $s(L^\infty(\mathcal G^{(0)}, \nu))$ is central in $L^\infty(\mathcal G, \mu)$). So, we get that $(b,\tilde{\ga})$ is outer. 
\end{proof}

\subsection{Theorem}
\label{thtensor}
{\it Let $\mathcal G$ be locally compact separable infinite groupoid; then $\mathcal G$ has an outer action on a hyperfinite semi-finite von Neumann algebra. }

\begin{proof}
Let us apply \ref{proptensor} to the action constructed in \ref{notgroupoid}. \end{proof}

\subsection{Theorem}
\label{finite}
{\it Let $\gG$ be a measured quantum groupoid, and let $(b, \ga)$ be an outer action of $\gG$ on a finite von Neumann algebra $A$; let $\theta$ be a faithful tracial state on $A$; then :
\newline
(i) there exists a positive non singular operator $h$ affiliated to the center of $N$, such that $(D\theta\circ b: D\nu)_t=h^{it}$; moreover, for all $x\in M$, we have $\tau_t(x)=\alpha(h^{-it})\beta(h^{it})x\alpha(h^{it})\beta(h^{-it})$. 
\newline
(ii) there exists a normal faithful conditional expectation from $A$ onto $b(N)$, which is invariant under $\ga$. 
}

\begin{proof}
As $\theta\circ b$ is a trace on $N$, we get that there exists a normal faithful conditional expectation $E$ from $A$ onto $b(N)$; therefore, the action $(b, \ga)$ is weighted, and we can apply \ref{semi-finite}, from which we get that $\nu$ is a trace, and, therfore, the existence of the operator $h$ defining the Radon-Nikodym derivative between $\theta\circ b$ and $\nu$. Moreover, as, for any positive $n$ in $N$, we have $\theta\circ b(n)=\nu(hn)$, we get that there exists also a normal semi-finite faithful operator-valued weight $\mathcal T$ from $A$ onto $b(N)$ such that $\theta=\nu\circ b^{-1}\circ \mathcal T$, which verify, for all positive $x$ in $A$, $\mathcal T(x)=b(h)E(x)$; moreover, using then \ref{rhoinv}, we get, for any $x\in\gM_{\mathcal T}^+$ :
\begin{eqnarray*}
(E\underset{N}{_b*_\alpha}id)\ga(x)(1\underset{N}{_b\otimes_\alpha}\alpha(h))
&=&
(E\underset{N}{_b*_\alpha}id)\ga(x)(b(h)\underset{N}{_b\otimes_\alpha}1)\\
&=&
(\mathcal T\underset{N}{_b*_\alpha}id)\ga(x)\\
&=&
1\underset{N}{_b\otimes_\alpha}(\theta\underset{N}{_b*_\alpha}id)\ga(x)\\
&=&
1\underset{N}{_b\otimes_\alpha}\beta\circ b^{-1}(\mathcal T(x))\rho^{-1}\\
&=&
1\underset{N}{_b\otimes_\alpha}\beta\circ b^{-1}(E(x))\beta(h)\rho^{-1}
\end{eqnarray*}
And, making now $x$ increase towards $1$, we get that $\rho=\beta(h)\alpha(h^{-1})$, which is (i).
\newline
Using this result in the calculation above, we get, for any $x\in\gM_{\mathcal T}^+$ :
\[(E\underset{N}{_b*_\alpha}id)\ga(x)=1\underset{N}{_b\otimes_\alpha}\beta\circ b^{-1}(E(x))=\ga(E(x))\]
which, by increasing limits, remains true for any positive $x$ in $A$; which is (ii). 

\end{proof}

\subsection{Theorem}
\label{finite2}
{\it Let $\gG$ be a measured quantum groupoid, and let $(b, \ga)$ be an outer action on a von Neumann algebra $A$; let us suppose that the crossed-product $A\rtimes_\ga\gG$ is a finite von Neumann algebra. Then :
\newline
(i) $\widehat{\gG}$ is a measured quantum groupoid of compact type, in the sense of (\cite{L}, 13.2) and (\cite{E3}, 5.11), i.e. there exists a left-invariant normal conditional expectation on the Hopf-bimodule $(\widehat{M}, N, \alpha, \hat{\beta}, \widehat{\Gamma})$, which implies that $\hat{\delta}=\lambda=1$. Moreover, we get also that $M$ is semi-finite. 
\newline
(ii) we have :
\[Z(A)\underset{N}{_b\otimes_\alpha}\mathbb{C}=Z(A\rtimes_\ga\gG)\]
\[\alpha(N)\cap Z(M)=\alpha(N)\cap Z(\widehat{M})\]
Therefore, $A$ is a factor if and only if $A\rtimes_\ga\gG$ is a factor, and $\gG$ is connected if and only if $\widehat{\gG}$ is connected. 
\newline
(iii) let us suppose $A$ is a factor; then, $\gG$ is finite dimensional if and only if the depth 2 inclusion $\ga(A)\subset A\rtimes_\ga\gG$ is of finite index. }

\begin{proof}
Let $\theta$ be a faithful tracial normal state on $A\rtimes_\ga\gG$; then, the restriction of $\theta$ to $\ga(A)$ is also a faithful  normal state, and, there exists a normal faithful conditional expectation $E$ from $A\rtimes_\ga\gG$ onto $\ga(A)$. On the other hand, using (\cite{E5}, 9.8), we get that theere exists a normal semi-finite faithful operator-valued weight $T_{\tilde{\ga}}$ from $A\rtimes_\ga\gG$ onto $\ga(A)$, and we get that $(DT_{\tilde{\ga}}:DE)_t$ belongs to $A\rtimes_\ga\gG\cap\ga(A)'$, which is, as $\ga$ is outer, equal to $1\underset{N}{_b\otimes_\alpha}\hat{\alpha}(N)$. 
\newline
Moreover, let us consider $\theta\circ\ga$ which is is a faithful tracial normal state on $A$ (and allows us to apply \ref{finite}); we have :
\[(DT_{\tilde{\ga}}:DE)_t=(D\widetilde{\theta\circ\ga}:D\theta)_t\]
and, as $\theta$ is a trace, we get that there exists $k$ positive invertible affiliated to $N$ such that :
\[(DT_{\tilde{\ga}}:DE)_t=1\underset{N}{_b\otimes_\alpha}\hat{\alpha}(k^{it})\]
which gives that, for any positive $X$ in $A\rtimes_\ga\gG$, we have :
\[T_{\tilde{\ga}}(X)=E((1\underset{N}{_b\otimes_\alpha}\hat{\alpha}(k^{1/2}))X(1\underset{N}{_b\otimes_\alpha}\hat{\alpha}(k^{1/2})))\]
and, taking now a positive $y$ in $\widehat{M}'$, we have, using (\cite{E5}, 9.8) :
\[1\underset{N}{_b\otimes_\alpha}\hat{T}^c(y)=E(1\underset{N}{_b\otimes_\alpha}\hat{\alpha}(k^{1/2})y\hat{\alpha}(k^{1/2}))\]
From which, by taking the restriction of $E$ to $1\underset{N}{_b\otimes_\alpha}\widehat{M}'$, we get a normal faithful conditional expectation $F$ from $\widehat{M}'$ onto $\beta(N)$ which satisfies, for all positive $y$ in $\widehat{M}'$ :
\[\hat{T}^c(y)=F(\hat{\alpha}(k^{1/2})y\hat{\alpha}(k^{1/2}))\]
Or, equivalently, we get the existence of a normal faithful conditional expectation $G$ from $\widehat{M}$ onto $\alpha(N)$ such that, for all positive $z$ in $\widehat{M}$, we have :
\[\hat{T}(z)=G(\hat{\beta}(k^{1/2})z\hat{\beta}(k^{1/2}))\]
It is then straightforward to verify that $G$ is left-invariant, which gives the beginning of (i). As $\widehat{\gG}$ is of compact type, we have $\widehat{\Phi}=\widehat{\Phi}\circ\hat{R}$, which implies $\hat{\delta}=\hat{\lambda}=1$; as  $\hat{\lambda}=\lambda^{-1}$, (\cite{E5}, 3.10 (vii)), we get that $\lambda=1$. Moreover, we get also that $\Delta_\Phi=P$, as $P^{it}$ is the standard implementation of $\tau_t$, which is, thanks to \ref{finite}, equal to the interior automorphism implemented by $\alpha(h^{-it})\beta(h^{it})$, where $h$ is defined as $\theta(1\underset{N}{_b\otimes_\alpha}\beta(n))=\nu(hn)$, for all positive $n$ in $N$. So, $\sigma_t^\Phi=\tau_t$ is interior, which finishes the proof of (i). 
\newline
By applying again \ref{finite}, we get that the action $\ga$ is weighted, so that we may apply \ref{semi-finite} to the action $\ga$; but, we may also apply \ref{finite} and \ref{semi-finite} to the action $\tilde{\ga}$; so, we get that :
\[Z(A)\underset{N}{_b\otimes_\alpha}\mathbb{C}\subset Z(A\rtimes_\ga\gG)\subset Z(A\underset{N}{_b*_\alpha}\mathcal L(H))=Z(A)\underset{N}{_b\otimes_\alpha}\mathbb{C}\]
and that :
\[\alpha(N)\cap Z(M)\subset \alpha(N)\cap Z(\widehat{M})\subset \alpha(N)\cap Z(M)\]
which finishes the proof of (ii). 
\newline
When $A$ is a factor, it is well known (\cite{GHJ}, 4.6.2) that, if the inclusion $\ga(A)\subset A\rtimes_\ga\gG$ is of finite index, then all relative commutants in the tower are finite-dimensional; in particular $A\underset{N}{_b*_\alpha}\mathcal L(H)\cap \ga(A)'$ is finite-dimensional, and we get that $\gG$ is finite-dimensional by \ref{thout}(ii). Conversely, if $\gG$ is finite-dimensional, so is $H$, and the factor $A\underset{N}{_b*_\alpha}\mathcal L(H)$ is finite, so is $\ga(A)'$, and the index of $\ga(A)\subset A\rtimes_\ga\gG$ is finite (\cite{J}, 2.1.7). 
\end{proof}

\subsection{Remarks}
\label{remarks}
(i) in (\cite{NSzW} 4.2.5) was proved that any connected finite dimensional quantum groupoid outerly acting on a factor is biconnected. 
\newline
(ii) in \cite{NV2} was proved that any depth 2 finite index subfactor of the hyperfinite $II_1$ factor $\mathcal R$ leads to a biconnected finite dimensional quantum groupoid outerly acting on $\mathcal R$ (such that the subfactor is the algebra of invariants elements by this action). 
\newline
(iii) in \cite{D} was proved that any finite dimensional biconnected quantum groupoid, whose antipode is involutive on the two copies of the basis (called target and source Cartan subalgebras, or target and source co-unital subalgebras) is outerly acting on $\mathcal R$; this hypothesis on the antipode is equivalent to the fact of having a finite normal quasi-invariant trace on the basis). 

\subsection{Examples}
\label{exs}
(i) as recalled in \ref{remarks}(iii), any connected measured quantum groupoid $\gG=(N, M, \alpha, \beta, \Gamma, T, T', \nu)$ such that $dim M<\infty$, and $\nu$ is a trace has an outer action $\ga$ on $\mathcal R$ such that the crossed-product $\mathcal R\rtimes_\ga\gG$ is isomorphic to $\mathcal R$, and the index of the inclusion $\ga(\mathcal R)\subset \mathcal R\rtimes_\ga\gG$ is finite. Moreover, by \ref{remarks}(ii), these are the only outer actions of finite index of measured quantum groupoids on $\mathcal R$.
\newline
(ii) let $\tau$ be the tracial normal faithful state on $\mathcal R$, and let us write $H$ for $H_\tau$, $J$ for $J_\tau$, $Tr$ the canonical semi-finite faithful trace on $\mathcal L(H)$. Let us denote by $T_\tau$ the normal faithful semi-finite operator-valued weight from $\mathcal L(H)$ onto $\mathcal R$, such that $\tau\circ T_\tau=Tr$. Let us recall the construction of the "$\mathcal R$-quantum groupoid", as made in (\cite{L},14). Let us consider the von Neumann algebra $\mathcal R^o\otimes \mathcal R$, equipped with its canonical structure of $\mathcal R$-bimodule, i.e. we define, for any $x\in\mathcal R$, $\alpha(x)=1\otimes x$, and $\beta(x)=x^o\otimes 1$; then, the fiber product $(\mathcal R^o\otimes \mathcal R)\underset{\mathcal R}{_\beta*_\alpha}(\mathcal R^o\otimes \mathcal R)$ is canonically isomorphic to $\mathcal R^o\otimes \mathcal R$; moreover $\gG(\mathcal R)=(\mathcal R, \mathcal R^o\otimes \mathcal R, \alpha, \beta, id, id\otimes\tau, \tau^o\otimes id, \tau)$ is a measured quantum groupoid, which is of compact type, because $id\otimes\tau$ is a conditional expectation. Constructing its dual, we obtain the von Neumann algebra $\mathcal L(H)$, with its canonical structure of $\mathcal R$-bimodule ; let us write $id_{|\mathcal R}$ for the inclusion of $\mathcal R$ into $\mathcal L(H)$, and $id^o_{|\mathcal R}$ for the canonical anti-homomorphism $x\mapsto Jx^*J$ from $\mathcal R$ into $\mathcal L(H)$; then, we get that the fiber product $\mathcal L(H)\underset{\mathcal R}{_{id^o_{|\mathcal R}}*_{id_{|\mathcal R}}}\mathcal L(H)$ is canonically isomorphic to $\mathcal L(H)$, and we get that $\widehat{\gG(\mathcal R)}$ is equal to $(\mathcal R, \mathcal L(H), id_{|\mathcal R}, id^o_{|\mathcal R}, id, T_\tau, (T_\tau)^o, \tau)$, where $id$ means the identity of $\mathcal L(H)$, and $(T_\tau)^o$ is defined, for any $y\in\mathcal L(H)$, by $(T_\tau)^o(y)=JT_\tau(JyJ)J$. 
\newline
Let us consider now the trivial action $(id^o_{|\mathcal R}, id_{|\mathcal R})$ of $\widehat{\gG(\mathcal R)}$ on $\mathcal R$ (\cite{E5}, 6.2): it's crossed-product (\cite{E5}, 9.5) is the von Neumann algebra of $\gG(\mathcal R)^c$, i.e. $\mathcal R\otimes \mathcal R^o$, which is isomorphic to $\mathcal R$; moreover, the relative commutant of $\mathcal R$ into the crossed-product is $\mathcal R^o$, and, so, this action is outer. The inclusion of $\mathcal R$ into the crossed-product is isomorphic to $\mathcal R\otimes \mathbb{C}\subset \mathcal R\otimes\mathcal R^o$, which is of infinite index (\cite{J}, 2.1.19). 
\newline
(iii) Let $\gG$ be a finite dimensional measured quantum groupoid, with a relatively invariant trace $\nu$ on the basis $N$, and $(b, \ga)$ an outer action of $\gG$ on $\mathcal R$; let $\widehat{\gG(\mathcal R)}$ be the dual $\mathcal R$-quantum groupoid, and $(id, id)$ its trivial action on $\mathcal R$, which is outer by (ii). We can construct now the measured quantum groupoid $\gG\oplus\widehat{\gG(\mathcal R)}$ (\ref{exMQG}(v)) and construct the action $(b\oplus id, \ga\oplus id)$ of $\gG\oplus\widehat{\gG(\mathcal R)}$ on $\mathcal R\oplus\mathcal R$. We obtain this way (\ref{lcqg}(ii)) an outer action of $\gG\oplus\widehat{\gG(\mathcal R)}$ on $\mathcal R\oplus\mathcal R$ (or, equivalently, on $\mathcal R$), whose crossed-product is also isomorphic to $\mathcal R$. This action is clearly of infinite index. 
\newline
(iv) Let $\bf{G}$ be a locally compact quantum group having a strictly outer action (in the sense of \cite{V2}) on $\mathcal R$; for instance, any locally compact group $G$ (\cite{V2}, 5.2), or any amenable Kac algebra of discrete type (\cite{V2}, 8.1) (or, equivalently, using (\cite{To}, 3.17), any Kac algebra of discrete type, such that the underlying von Neumann algebra of the dual Kac algebra of compact type is injective), or, by duality, any Kac algebra of compact type whose undelying von Neumann algebra is injective). Using again \ref{lcqg}(ii) and the example given in (iii), we obtain the existence of an outer action of $\bf{G}\oplus\gG\oplus\widehat{\gG(\mathcal R)}$ on $\mathcal R$, for any such $\bf{G}$, any finite dimensional measured quantum groupoid $\gG$, with a relative invariant trace on the basis. Using \ref{finite2}(i), we get that its crossed-product is finite if and only if $\bf{G}$ is of discrete type. This crossed-product is then a finite factor, by \ref{finite2}(ii), and, in the case when this Kac algebra of discrete type has an invariant mean, (which, by (\cite{To}, 3.17), implies that the underlying von Neumann algebra of the dual Kac algebra of compact type is injective), we get that the crossed-product is injective and is therefore isomorphic to $\mathcal R$, and this inclusion will be of infinite index. 
\newline
(v) Let's now give examples of outer actions of measured quantum groupoids on semi-finite von Neumann algebras. We had got in \ref{thtensor} that any locally compact separable infinite groupoid has an outer action on a hyperfinite semi-finite von Neumann algebra; therefore, using (iv), we easily get that $\gG(\mathcal G)\oplus\bf{G}\oplus\gG\oplus\widehat{\gG(\mathcal R)}$ has on outer action on a hyperfinite semi-finite von Neumann algebra, where $\gG(\mathcal G)$ is the measured quantum groupoid constructed from $\mathcal G$, $\bf{G}$ is any locally compact quantum group having a strictly outer action on $\mathcal R$, $\gG$ is any finite dimensional measured quantum groupoid, with a relatively invariant trace on the basis, and $\widehat{\gG(\mathcal R)}$ has been defined in (ii). 


\end{document}